\documentclass[12pt]{article}
\usepackage{amsmath}
\usepackage{amssymb}
\usepackage{amscd}
\newtheorem{theoreme}{Th\'eor\`eme}[section]
\newtheorem{proposition}[theoreme]{Proposition}
\newtheorem{corollaire}[theoreme]{Corollaire}
\newtheorem{lemme}[theoreme]{Lemme}
\newtheorem{exemple}[theoreme]{Exemple}
\newtheorem{remarque}[theoreme]{Remarque}
\newenvironment{preuve}{\begin{trivlist} \item[]{\it Preuve---}}
{\par\hfill $\square$\end{trivlist}}
\renewcommand{\P}{\mathbb{P}}
\newcommand{\C}{\mathbb{C}}
\newcommand{\R}{\mathbb{R}}
\newcommand{\N}{\mathbb{N}}
\newcommand{\Z}{\mathbb{Z}}
\newcommand{\T}{{\rm T}}
\newcommand{\id}{{\rm id}}
\newcommand{\M}{{\cal M}}
\newcommand{\TF}{{\ \!\! ^\tau\!}}
\renewcommand{\H}{{\cal H}}
\newcommand{\G}{{\cal G}}
\newcommand{\A}{{\cal A}}

\newcommand{\E}{{\cal E}}
\newcommand{\D}{{\cal D}}

\newcommand{\mult}{{\rm mult}}
\newcommand{\const}{{\rm const}}

\renewcommand{\Re}{{\rm Re}}
\renewcommand{\Im}{{\rm Im}}
\renewcommand{\O}{{\rm O}}
\renewcommand{\o}{{\rm o}}
\newcommand{\ahead}{\par \hspace{0.75cm}}
\title{Sur les endomorphismes holomorphes permutables de
$\P^k$}
\author{Tien-Cuong Dinh et Nessim Sibony}
\begin{document}
\maketitle
\begin{abstract} Let $f_1$, $f_2$ be holomorphic endomorphisms of
  $\P^k$, of degrees $d_1\geq 2$, $d_2\geq 2$. Assume that $f_1\circ
  f_2=f_2\circ f_1$ and that $d_1^{n_1}\not =d_2^{n_2}$ for all
  integers $n_1$, $n_2$. We then show that $f_j$ are critically
  finite. Moreover there is an orbifold $(\P^k,n)$ such that $f_1$,
  $f_2$ are coverings of $(\P^k,n)$. In the $\P^2$ case we give the
  list of commuting pairs satisfying the above conditions.
\end{abstract}
\section{Introduction}
\ahead
Dans \cite{Smale}, S. Smale pose le probl\`eme suivant. Etant donn\'e
une vari\'et\'e compacte $M$, est-ce que tout diff\'eomorphisme de $M$
peut \^etre approch\'e par des diff\'eomorphismes qui commutent
seulement avec leurs it\'er\'es? Il ajoute ``I find this problem
interesting in that it gives some focus in the dark realm, beyond
hyperbolicity, where even the problems are hard to pose clearly''.
\par
Dans le cadre de la dynamique des applications rationnelles de $\P^1$
le probl\`eme a vivement int\'eress\'e Fatou \cite{Fatou} et Julia
\cite{Julia}. De fait ils ont r\'esolu l'\'equation fonctionnelle: 
\begin{equation}
f_1\circ f_2  =  f_2\circ f_1 
\end{equation}
o\`u $f_1$ et $f_2$ sont des polyn\^omes d'une variable
complexe de degr\'e $d_1\geq 2$ et $d_2\geq 2$
\cite{Fatou,Julia}. Ritt \cite{Ritt} a r\'esolu cette \'equation
fonctionnelle pour les fractions rationnelles. Fatou-Julia
consid\'eraient que l'un des buts de la th\'eorie de l'it\'eration
\'etait l'investigation des \'equations fonctionnelles.
\par
Lorsque deux applications rationnelles sont {\it permutables} 
(c'est-\`a-dire si 
elles v\'erifient la relation (1)), les objets dynamiques qui
leur sont associ\'es sont fortement li\'es. En particulier, elles ont
m\^eme ensemble de Julia et m\^eme mesure d'\'equilibre. Utilisant ces
notions dynamiques Eremenko \cite{Eremenko} a donn\'e une version
nouvelle des r\'esultats de Ritt. Il utilise en particulier la notion
d'orbifold consid\'er\'ee par Thurston \cite{Thurston}.
\par
Dans cet article nous nous int\'eressons \`a l'\'equation (1) lorsque
$f_1$, $f_2$ sont des endomorphismes holomorphes de $\P^k$.
\par
Pour cela nous nous utilisons sur les progr\`es r\'ecents de la
th\'eorie de l'it\'eration des applications rationnelles de
$\P^k$. Nous renvoyons \`a 
\cite{BedfordSmillie, Fornaess, Sibony} pour une description de ces
d\'eveloppements. L'expos\'e \cite{Sibony} est adapt\'e \`a nos besoins.
\par
Si on pose $h^m:=h\circ\cdots\circ h$ ($m$ fois) et
$h^n:=h\circ\cdots\circ h$ ($n$ fois). 
Il est clair que $h^m$ et $h^n$ sont permutables lorsque $h$ est un
endomorphisme holomorphe de $\P^k$. Nous limiterons notre \'etude au
cas o\`u
\begin{equation}
f_1^{n_1}\not=f_2^{n_2} \mbox{ pour tous nombres naturels non nuls } n_1
\mbox{ et } n_2
\end{equation}
Cette condition est vraie en particulier si
\begin{equation}
d_1^{n_1}\not = d_2^{n_2} \mbox{ pour tous nombres naturels non nuls }
n_1 \mbox{ et } n_2
\end{equation}
Notons $w:=[w_0:\cdots:w_k]$ (resp. $z:=(z_1,\ldots,z_k)$) les
coordonn\'ees homog\`enes (resp. les coordonn\'ees affines) de
$\P^k$ o\`u $z_s:=w_s/w_0$.
\par
Dans le cas de dimension 1 ($k=1$), les solutions de (1)(2)
sont (pour une certaine coordonn\'ee $z$, \cite{Eremenko}):
\begin{enumerate}
\item $f_1(z)=z^{\pm d_1}$ et $f_2(z)=\lambda z^{\pm d_2}$ avec
$\lambda\not = 0$ convenablement choisi;
\item $f_1(z)=\pm\T_{d_1}(z)$ et $f_2(z)=\pm \T_{d_2}(z)$ avec
les signes $\pm$ convenables o\`u $\T_{d_i}(\cos z):=\cos(d_iz)$
est le {\it polyn\^ome de Tchebychev} de degr\'e $d_i$;
\item $f_1$ et $f_2$ sont des {\it applications de Latt\`es}, c'est-\`a-dire il
existe une application holomorphe surjective
$\varphi:\ \C\longrightarrow
\P^1$, des applications affines $\Lambda_i$ et un groupe
d'automorphismes affines discret $\A$ de $\C$, agissant
transitivement sur les fibres de $\varphi$ tels que
$f_i\circ\varphi=\varphi\circ\Lambda_i$. 
\end{enumerate}
\par
Notre r\'esultat principal est le th\'eor\`eme suivant, qui dans
le cas de dimension 1, permet de retrouver les trois solutions d\'ecrites
ci-dessus:
\begin{theoreme} Soient $f_1$ et $f_2$ deux endomorphismes
holomorphes permutables de degr\'es $d_1\geq 2$ et $d_2\geq 2$
de $\P^k$. Supposons
que $d_1^{n_1}\not=d_2^{n_2}$ pour tous $n_1$, $n_2$ entiers strictement
positifs. Alors il
existe une application holomorphe $\varphi: \C^k\longrightarrow
\P^k$ et des applications affines holomorphes
$\Lambda_1$, $\Lambda_2$ 
tels que $\varphi(\C^k)$ soit un ouvert de compl\'ement pluripolaire 
de $\P^k$ et tels que
$f_i\circ\varphi=\varphi\circ\Lambda_i$. De plus, il existe  
un groupe discret d'applications affines holomorphes
$\A$ de $\C^k$
agissant transitivement
sur les fibres de $\varphi$.
\end{theoreme}
\begin{remarque} \rm
Les applications $\varphi$ qui interviennent dans ce
  th\'eor\`eme sont donc invariantes par le groupe $\A$. 
\end{remarque}
\par
Soit $X$ une vari\'et\'e complexe. Notons $\H(X)$
l'espace des sous-ensembles analytiques
irr\'eductibles de codimension 1 de $X$.
On appelle {\it orbifold} un
couple $(X,n)$ o\`u $n$ est une
fonction d\'efinie sur $\H(X)$ \`a valeurs dans
$\N^+\cup\{\infty\}$ \'egale \`a 1 sauf sur une
famille localement finie de sous-ensembles analytiques de $X$.
{\it Un
rev\^etement d'orbifolds} $\pi:\ (X_1,n_1)\longrightarrow
(X_2,n_2)$ est un rev\^etement ramifi\'e de $X_1\setminus
\bigcup_{n_1(H)=\infty}H$ dans $X_2\setminus\bigcup_{n_2(H)=
\infty}H$
v\'erifiant $\mult(\pi,H).n_1(H)=n_2(\pi(H))$ pour tout
$H\in\H(X_1)$, o\`u $\mult(\pi,H)$ d\'esigne la multiplicit\'e de $\pi$
en un point g\'en\'erique de $H$.
\par
Soit $f$ une application holomorphe de $\P^k$ dans
$\P^k$ d\'efinissant un rev\^etement d'un orbifold ${\cal
O}=(\P^k,n)$ dans lui-m\^eme. Alors, l'ensemble critique de $f$
est {\it pr\'ep\'eriodique}, c'est-\`a-dire $f^n({\cal C}_f)=f^m({\cal C}_f)$
pour certains $0\leq n<m$ o\`u ${\cal C}_f$ d\'esigne l'ensemble
critique de $f$. On dit qu'une telle application est {\it
critiquement finie}. 
Dans le cas de dimension 1, le corollaire suivant
se r\'eduit au th\'eor\`eme de Fatou-Julia-Ritt \cite{Eremenko}:
\begin{corollaire} Sous l'hypoth\`ese du th\'eor\`eme
1.1, il existe un orbifold ${\cal O}=(\P^k,n)$ tel que les
applications $f_1$ et $f_2$ d\'efinissent des rev\^etements de
${\cal O}$ dans lui-m\^eme. En particulier, $f_1$ et $f_2$ sont
critiquement finies. 
\end{corollaire}
Pour tout $k\geq 3$, 
le th\'eor\`eme 1.1 et le corollaire 1.3 ne sont plus vrais si l'on
remplace la condition (3) par la condition (2). Donnons un exemple.
\begin{exemple} \rm
Soient $h_1$ et $h_2$ deux fractions rationnelles de Latt\`es
permutables, de m\^eme degr\'e $d\geq 2$ v\'erifiant
$h_1^n\not=h_2^n$ pour tout $n\geq 1$. On peut prendre par exemple
$\varphi:\C\longrightarrow \P^1$ la fonction elliptique de Weierstrass
de p\'eriodes 1 et $i$; $\A$ le groupe engendr\'e par les
automorphismes $z\mapsto z+1$, $z\mapsto z+i$ et $z\mapsto -z$;
$\Lambda_1(z):= (1+2i)z$; $\Lambda_2(z):=(1-2i)z$; $d=5$ et les
applications $h_1$, $h_2$ satisfaisant $h_1\circ\varphi=\varphi\circ
\Lambda_1$, $h_2\circ\varphi=\varphi\circ
\Lambda_2$ ({\it voir} \cite{Eremenko}).
Dans les coordonn\'ees
homog\`enes de $\P^1$,
on peut \'ecrire $h_i=[P_i:Q_i]$ o\`u $P_i$ et $Q_i$
sont des polyn\^omes homog\`enes de deux variables de degr\'e
$d$. Consid\'erons les deux endomorphismes holomorphes de $\C^3$
d\'efinis par $f_1(z)=(P_1(z_1,z_2),Q_1(z_1,z_2),R(z_3))$ et
$f_2(z)=(\lambda P_2(z_1,z_2),\lambda Q_2(z_1,z_2),R(z_3))$
o\`u $R$ est un polyn\^ome de degr\'e $d$ et $\lambda$ est une
constante non nulle. Ces endomorphismes se prolongent en des
endomorphismes holomorphes de $\P^3$ v\'erifiant (2).
Pour des $\lambda$ convenables, ces endomorphismes sont permutables.
Si $R$ n'est pas critiquement fini, $f_1$ et $f_2$ ne sont pas
critiquement finis. Le Th\'eor\`eme 1.1 n'est pas valide dans ce cas
car le groupe $\A$ n'existe pas.
\end{exemple}
\par
Les exemples suivants sont des solutions du probl\`eme
(1)(2). Pour simplifier les notations, nous nous limitons au cas de $\P^2$.
\begin{exemple} \rm
Soient $h_1$ et $h_2$ des
endomorphismes holomorphes permutables de $\P^1$. Dans les
coordonn\'ees homog\`enes de $\P^1$, il existe des polyn\^omes
homog\`enes \`a deux variables $P_i$ et $Q_i$ tels que
$h_i=[P_i:Q_i]$ pour $i=1$ ou 2. Les
endomorphismes holomorphes $f_i$
de $\P^2$ d\'efinis en coordonn\'ees affines $z$ par
$f_1(z):=(P_1(z_1,z_2),Q_1(z_1,z_2))$ et $f_2(z):=(\lambda
P_2(z_1,z_2),\lambda Q_2(z_1,z_2))$ sont permutables pour des
constantes $\lambda\not=0$ convenables.
\end{exemple}
\begin{exemple} \rm
Les endomorphismes holomorphes
$f_1(w):=[w_{\alpha_0}^{d_1}:
w_{\alpha_1}^{d_1}:w_{\alpha_2}^{d_1}]$ et
$f_2(w):=[\lambda_0 w_{\nu_0}^{d_1}:
\lambda_1 w_{\nu_1}^{d_1}:\lambda_2 w_{\nu_2}^{d_1}]$ sont
permutables lorsque $\{\alpha_0,\alpha_1,\alpha_2\}$ et
$\{\nu_0,\nu_1,\nu_2\}$ sont des permutations de $\{0,1,2\}$ et
$\lambda_0$, $\lambda_1$, $\lambda_2$ sont des 
constantes non nulles convenablement choisies.
\end{exemple}
\begin{exemple} \rm
Les endomorphismes $f_1(z_1,z_2):=(z_1^{\pm d_1},
  \pm \T_{d_1}(z_2))$ et $f_2(z_1,z_2):=(\lambda z_1^{\pm d_2},
  \pm \T_{d_2}(z_2))$ 
sont permutables lorsque la constante $\lambda\not =0$ 
 et les signes $\pm$ sont convenablement choisis.
\end{exemple}
\begin{exemple} \rm
Les endomorphismes $f_i(z_1,z_2):=(\pm
  \T_{d_i}(z_1), \pm \T_{d_i} (z_2))$ ou\break
$(\pm \T_{d_i}(z_2), \pm \T_{d_i}(z_1))$
  sont permutables lorsque les signes $\pm$ sont convenablement
choisis.
\end{exemple}
\begin{exemple} \rm
Soient $h_1$, $h_2$ deux endomorphismes
holomorphes permutables 
  de $\mathbb{P}^1$. Soit $\pi$ l'application holomorphe de
  $\mathbb{P}^1\times \mathbb{P}^1$ dans $\mathbb{P}^2$ qui d\'efinit
  un rev\^etement ramifi\'e \`a deux feuillets tel que
  $\pi(x,y)=\pi(y,x)$. 
Posons $F_i(a,b):=(h_i(a),h_i(b))$ deux
  endomorphismes de $\mathbb{P}^1\times\mathbb{P}^1$. Il existe des
  endomorphismes holomorphes permutables
$f_i$ de $\mathbb{P}^2$ tels que
  $f_i\circ\pi=\pi\circ F_i$. 
Lorsque les
  $h_i$ sont des polyn\^omes, $f_1$, $f_2$ sont polynomiaux.
Lorsque les $h_i$ sont des applications des Latt\`es, on obtient un cas
particulier de l'exemple 1.10.
\end{exemple}
\begin{exemple} \rm
Dans le contexte du th\'eor\`eme 1.1, lorsque
$\A$ est un groupe cristallographique complexe (c'est-\`a-dire un
groupe co-compact d'isom\'etries complexes de $\C^k$), on dite que $f_1$ et
$f_2$ sont des {\it applications 
de Latt\`es g\'en\'eralis\'ees}.
Lorsque $\Lambda_1\circ\Lambda_2=\Lambda_2\circ\Lambda_1$ modulo
$\A$, les applications $f_1$ et $f_2$ sont permutables. 
Certains sous-familles d'applications de Latt\`es
g\'en\'eralis\'ees  sont
pr\'ecisement d\'ecrites dans \cite{BertelootLoeb}. 
\end{exemple}
\par
Nous obtenons le corollaire suivant qui
g\'en\'eralise le r\'esultat principal de \cite{Dinh2}:
\begin{corollaire} Sous l'hypoth\`ese du th\'eor\`eme
1.1, si $k=2$, le couple $(f_1,f_2)$ est \'egal,
dans un syst\`eme de coordonn\'ees
convenable de $\P^2$, \`a l'un des couples d\'ecrits dans les
exemples pr\'ec\'edents.
\end{corollaire}
\begin{remarque} \rm
Le m\^eme probl\`eme pour les automorphismes polynomiaux de
$\C^2$ est r\'esolu par Lamy \cite{Lamy}. 
Veselov \cite{Veselov} a donn\'e aussi 
une liste d'applications polynomiales
permutables dites {\it applications de Tchebychev}.
\end{remarque}
\par
Notre approche est semblable \`a celle adopt\'ee par Fatou, Julia et
Eremenko \`a une variable.
\par
Soient $f_1$, $f_2$ deux endomorphismes permutables de $\P^k$ de
degr\'e alg\'ebrique respectifs $d_1$ et $d_2$. On montre \`a l'aide
d'un r\'esultat de Briend-Duval \cite{BriendDuval1} qu'ils ont une
infinit\'e de points p\'eriodiques r\'epulsifs communs. Supposons pour
simplifier que $f_1(a)=f_2(a)=a$ soit un tel point. On montre qu'il
existe une application (de Poincar\'e) $\varphi:\C^k\longrightarrow
\P^k$ telle que $\varphi(0)=a$, $\varphi'(0)$ inversible et 
$$f_i\circ \varphi=\varphi\circ\Lambda_i$$
o\`u les $\Lambda_i$ sont des applications triangulaires. Le
probl\`eme est d'analyser les applications $\Lambda_i$. Pour cela nous
utilisons syst\'ematiquement la fonction de Green commune associ\'ee
aux deux applications $f_i$. On montre que sous l'hypoth\`ese (3) le
groupe de Lie $\Gamma$ engendr\'e par $\Lambda_1$, $\Lambda_2$ contient un
sous-groupe isomorphe \`a $\R$. Sous l'hypoth\`ese (2), si
$d_1^{n_1}=d_2^{n_2}$ pour certains entiers strictement positifs
$n_1$, $n_2$ alors ce groupe contient un sous-groupe compact non
discret. Apr\`es avoir observ\'e que $f_1$ et $f_2$ ont m\^eme mesure
d'\'equilibre $\mu$, de support $J_k$, on montre sous l'hypoth\`ese
(3) que $J_k$ contient un ouvert qui est une vari\'et\'e 
r\'eelle analytique. En effet, $J_k$ est ``lamin\'e'' par les images
des orbites de $\Gamma$ (par l'application $\varphi$)
et on peut fabriquer des
``laminations'' diff\'erentes gr\^ace \`a des points p\'eriodiques
communs diff\'erents. 
On est donc amen\'e \`a \'etudier la structure
d'un endomorphisme $f$ de $\P^k$ pour lequel le support $J_k$ de la mesure
d'\'equilibre contient un ouvert $J_k\cap \Omega$ qui est 
une vari\'et\'e r\'eelle
analytique. C'est ce que nous faisons au paragraphe 5. 
\par
On consid\`ere une application de Poincar\'e $\varphi$ associ\'ee \`a
un point fixe r\'epulsif $b\in J_k\cap \Omega$ v\'erifiant $\varphi(0)=b$,
$\varphi'(0)$ inversible et satisfaisant l'\'equation
$$f\circ \varphi=\varphi\circ \Lambda.$$
On \'etudie alors les r\'esonnances possibles entre les valeurs
propres de $\Lambda'(0)$. Si $f$ est de degr\'e $d$, on montre que ces
valeurs propres sont \'egales \`a $\pm d$ ou de module $\sqrt{d}$. On
en d\'eduit que $J^*:=\varphi^{-1}(J_k)$ admet des \'equations de la forme 
$$\Im z''=q(z',\overline z') \mbox{ avec } z'\in\C^n, z''\in\C^{k-n}$$
o\`u $q$ est une application polynomiale homog\`ene de degr\'e deux. 
\par
On montre ensuite que si $g$ est un germe d'application holomorphe
v\'erifiant $\varphi\circ g=\varphi$ alors $g$ se prolonge en
application affine. Nous utilisons pour cela l'analyse de la fonction
$G^*:=G\circ \varphi$ o\`u $G$ est la fonction de Green associ\'ee \`a
$f$. L'une des difficult\'es pour montrer ce r\'esultat et construire
le groupe $\A$, est que $G^*$ n'est pas diff\'erentiable et qu'il faut
analyser les directions o\`u elle l'est dans un sens faible. Ces
outils mis en place on peut construire un orbifold associ\'e \`a
l'application $f$ comme au corollaire 1.3.
\par
On le voit que les endomorphismes de $\P^k$ de degr\'e $d_1$ qui sont
permutables \`a un endomorphisme de degr\'e $d_2$ avec $d_1^{n_1}\not
=d_2^{n_2}$ pour tout $(n_1,n_2)\not=(0,0)$ sont assez rigides et rares. 
\section{Germes d'applications holomorphes}
\ahead
Dans ce paragraphe, nous d\'emontrons que les germes
d'applications holomorphes permutables sont $\lambda$-triangulaires dans
un syst\`eme de coordonn\'ees convenables lorsque
$0$ est un point fixe r\'epulsif pour l'un d'eux.
Ce r\'esultat g\'en\'eralise un
r\'esultat similaire sur la triangulation des matrices permutables. 
\par
Notons $\G$ 
le semi-groupe des germes d'applications
holomorphes de $(\C^k,0)$ dans
$(\C^k,0)$.
Notons $\G^\times$ (resp. $\G^*$) 
le semi-groupe (resp. groupe) des germes d'applications
holomorphes \`a {\it fibres discr\`etes} (resp. inversibles) 
de $(\C^k,0)$ dans $(\C^k,0)$.
Soit $g\in\G$. On dit que $g$ est {\it triangulaire} si $g$ est inversible et
s'il est de la forme
$$g(z)=(\lambda_1 z_1,\lambda_2z_2+P_2(z),\ldots, \lambda_k
z_k+P_k(z))$$
o\`u  pour tout $2\leq j\leq k$ le polyn\^ome
$P_j$ est une combinaison lin\'eaire de mon\^omes
$z_1^{\alpha_1}\ldots z_{j-1}^{\alpha_{j-1}}$ dont le multi-indice
$(\alpha_1,\ldots,\alpha_{j-1})$ satisfait la relation
$\lambda_j=\lambda_1^{\alpha_1}
\ldots\lambda_{j-1}^{\alpha_{j-1}}$.
On dit que $0$ est un point {\it r\'epulsif} pour $g$ ou que $g$ est
{\it dilatant}
si toute valeur propre de $g'(0)$ est de module sup\'erieur \`a
1. D'apr\`es le th\'eor\`eme de Sternberg
\cite{Sternberg}, pour tout $g\in{\cal G}$ dilatant, il
existe $\varphi\in \G^*$ tel que
$\varphi^{-1}\circ g\circ\varphi$ soit triangulaire.
Si $g$ est triangulaire et dilatant
nous pouvons supposer que $1<|\lambda_1|\leq
|\lambda_2|\leq \cdots\leq |\lambda_k|$.
\par
Soit $\lambda:=(\lambda_1,\ldots,\lambda_k)\in(\C^*)^k$. Un
polyn\^ome $P(z)$ est dit {\it $\lambda$-homog\`ene} d'ordre $j$ s'il est
combinaison lin\'eaire de mon\^omes $z_1^{\alpha_1}\ldots
z_k^{\alpha_k}$ v\'erifiant $\lambda_j=
\lambda_1^{\alpha_1}\ldots\lambda_k^{\alpha_k}$. Il s'agit donc des
polyn\^omes qui v\'erifient $P(\lambda z)=\lambda_j P(z)$ o\`u
$\lambda z:=(\lambda_1 z_1,\ldots, \lambda_k z_k)$.
Un \'el\'ement $h=(h_1,\ldots,h_k)\in\G$
est dit {\it $\lambda$-triangulaire} si $h_j$ est
$\lambda$-homog\`ene d'ordre $j$ pour tout $1\leq j\leq k$. Autrement
dit les mon\^omes qui interviennent dans l'\'ecriture de $h$ sont ceux
donn\'es par les r\'esonnances de $\lambda$.
Notons
$\G_\lambda$ (resp. $\G_\lambda^*)$ l'ensemble des \'el\'ements
$\lambda$-triangulaires de $\G$ (resp. de $\G^*$).
On v\'erifie sans peine que si $|\lambda_j|>1$ pour tout
$1\leq j\leq k$, le semi-groupe $\G_\lambda$ est un espace complexe de
dimension finie et le groupe $\G_\lambda^*$ 
admet une structure naturelle de groupe de Lie. En
effet il y a seulement un nombre fini de r\'esonnances. 
\begin{proposition}
Soient $g_1$ et $g_2$ deux \'el\'ements permutables de ${\cal
G}^\times$. Supposons que $0$ soit r\'epulsif pour $g_1$.
Alors $g_2$ est inversible.
\end{proposition}
\begin{preuve}
Quitte \`a remplacer $g_2$ par $g_1^m\circ g_2$ pour $m$
suffisamment grand, on peut supposer
que les valeurs propres non nulles de $g_2'(0)$ sont de module
plus grand que $1$. Ceci est bien clair dans un syst\`eme de
coordonn\'ees o\`u les matrices permutables $g_1'(0)$ et
$g_2'(0)$ sont triangulaires. 
Supposons que $g_2'(0)$ ne soit pas inversible.
Soit $V$ la vari\'et\'e stable de $g_2$, c'est-\`a-dire $V:=\{z|\
\lim_{n\rightarrow \infty} g_2^n(z)=0\}$ \cite[p.27]{Ruelle}. Comme $g_1\circ
g_2=g_2\circ g_1$, $g_1(V)$ est \'egalement stable par $g_2$.
D'o\`u $V=g_1(V)$.
Quitte \`a remplacer $g_i$ par $g_i|_V$, on peut
supposer que les valeurs propres de $g_2'(0)$ sont nulles. En
rempla\c cant $g_2$ par $g_2^m$ avec $m$ suffisamment grand, on
peut supposer que $g_2'(0)=0$. D'apr\`es le th\'eor\`eme de Sternberg,
on peut supposer que $g_1$ est triangulaire. Comme $g_2$ est \`a
fibres discr\`etes,
$g_2(0,\ldots,0,z_k)$ n'est pas constant. La relation
$g_2=g_1^{-1}\circ g_2\circ g_1$ entra\^{\i}ne $g_2(0,\ldots
,0,z_k)=g_1^{-1}\circ g_2(0,\ldots,0,\lambda_kz_k)$ o\`u
les $\lambda_j$ sont les \'el\'ements de la diagonale principale de
$g_1'(0)$. On peut supposer que $1\leq |\lambda_1|\leq \cdots \leq
|\lambda_k|$. 
Notons $h_j$ les fonctions coordonn\'ees de
$g_2(0,\ldots,0,z_k)$. Soit $h_s$ la premi\`ere fonction non
nulle. Alors $h_s(z_k)=\lambda_s^{-1}h_s(\lambda_k z_k)$.
Si l'on consid\`ere la s\'erie de Taylor de $h_s$, la relation
pr\'ec\'edente contredit le fait que $h_s\not=0$, $h_s(0)=h_s'(0)=0$
et $1<|\lambda_s|\leq |\lambda_k|$.
\end{preuve}
\begin{proposition} Soient $g_1$ et $g_2$ deux \'el\'ements
permutables de $\G^\times$. 
Soient $\lambda_1,\ldots,\lambda_k$ les valeurs propres de
$g_1'(0)$ rang\'es dans l'ordre de croissance de leurs
modules. Posons $\lambda:=(\lambda_1,\ldots,\lambda_k)$.
Supposons que $0$ soit r\'epulsif
pour $g_1$. 
Alors il existe
$\varphi\in\G^*$ (dite application de
Poincar\'e) tel que 
$\varphi^{-1}\circ g_i\circ \varphi$ appartienne \`a $\G_\lambda^*$
pour $i=1$ ou $2$.
\end{proposition}
\begin{preuve} D'apr\`es la proposition 2.1, $g_2$ est inversible.
Comme $g_1\circ g_2=g_2\circ g_1$, les matrices $g_1'(0)$ et $g_2'(0)$
commutent. 
Quitte \`a faire un changement lin\'eaire de coordonn\'ees, on peut
supposer que les matrices $g_i'(0)$ sont triangulaires. De plus,
on peut supposer que les \'el\'ements de la diagonale principale
de $g_1'(0)$  sont $\lambda_1$, ..., $\lambda_k$. 
D'apr\`es le th\'eor\`eme de Sternberg, 
on peut supposer que
$g_1$ est $\lambda$-triangulaire. Il
suffit maintenant de 
montrer que $g_2$ l'est aussi.
On peut \'ecrire 
$g_1(z)=(\lambda_1 z_1,\lambda_2 z_2+P_2,\ldots,\lambda_k z_k+P_k)$
o\`u $P_j$ est un polyn\^ome $\lambda$-homog\`ene d'ordre $j$
en $z_1$, ..., $z_{j-1}$.
Posons $g_2=(h_1,\ldots, h_k)$. Soit
$s$ l'entier maximal tel que $h_j$ soit $\lambda$-homog\`ene d'ordre $j$
pour tout $1\leq j\leq s-1$. Montrons que $s=k+1$.
Supposons que $s\leq k$. La relation $g_1\circ g_2=g_2\circ g_1$
entra\^{\i}ne
$$h_s(\lambda_1 z_1,\ldots,\lambda_k z_k+P_k)=\lambda_s h_s+
P_s(h_1,\ldots ,h_{s-1}).$$
Soit $h$ la somme des termes dans le d\'eveloppement de Taylor de
$h_s$ qui ne sont pas $\lambda$-homog\`enes d'ordre $s$. 
Le choix de $s$ entra\^{\i}ne que $h\not
=0$ et que $P_s(h_1,\ldots,h_{s-1})$ est $\lambda$-homog\`ene d'ordre $s$.
On d\'eduit de l'\'equation
ci-dessus que 
$h\circ g_1=\lambda_sh$. Posons $\Lambda(z):=(\lambda_1
z_1,\ldots,\lambda_k z_k)$. Comme $g_1\circ\Lambda=\Lambda\circ g_1$,
on a $h\circ\Lambda^n\circ g_1=\lambda_s h\circ\Lambda^n$ pour tout
entier relatif $n$.
On choisit un mon\^ome $z_1^{n_1}\ldots z_k^{n_k}$ dont le
coefficient dans $h$ soit non nul et tel que 
$c:=|\lambda_1^{n_1}\ldots\lambda_k^{n_k}|$ soit minimal pour cette
propri\'et\'e. Il existe alors
une suite croissante d'entiers positifs $\{k_j\}$ telle que 
$c^{k_j} h\circ\Lambda^{-k_j}$ tende vers un polyn\^ome non nul
$P$. On a donc $P\circ g_1=\lambda_s P$.
Soit $z_1^{m_1}\ldots
z_k^{m_k}$ le terme dominant de $P$, c'est-\`a-dire pour tout terme
$z_1^{s_1}\ldots z_k^{s_k}$ de $P$
il existe $1\leq j\leq k$ v\'erifiant
$s_j<m_j$ et $s_i= m_i$ lorsque $i>j$. En identifiant les
coefficients de $z_1^{m_1}\ldots
z_k^{m_k}$ dans  $P\circ
g_1=\lambda_s P$ on obtient
$\lambda_1^{m_1}\ldots\lambda_k^{m_k}=\lambda_s$. Ceci contredit la
d\'efinition de $h$.
\end{preuve}
\par
Le lemme suivant sera utilis\'e pour prouver que $f_1$ et
$f_2$ poss\`edent une infinit\'e de points p\'eriodiques
r\'epulsifs communs.
\begin{lemme} Soient $f$, $g$ et $h$ trois
\'el\'ements de $\G^\times$ v\'erifiant
$f\circ h=h\circ g$. Si le point
$0$ est r\'epulsif pour $g$
alors il est r\'epulsif pour $f$.
\end{lemme}
\begin{preuve} Fixons un voisinage assez petit $V$ de $0$ et un
$n\geq 1$ tels  que $\overline V\subset g^n(V)$. Posons $U:=h(V)$.
Comme $h$ est \`a fibres discr\`etes, $U$ est un ouvert. La
relation $f^n\circ h=h\circ g^n$ entra\^{\i}ne $\overline
U\subset f^n(U)$. Cela entra\^{\i}ne que $0$ est r\'epulsif pour $f$. 
\end{preuve}
\section{Applications de Poincar\'e}
\ahead
Dans ce paragraphe, nous rappelons quelques outils de
la th\'eorie des syst\`emes dynamiques holomorphes: la fonction de
Green, les courants invariants, les ensembles de Julia. On pourra
trouver un expos\'e d\'etaill\'e dans \cite{Sibony}. Nous
allons d\'emontrer que $f_1$ et $f_2$ poss\`edent une infinit\'e
de points p\'eriodiques r\'epulsifs communs au voisinage desquels
des it\'er\'es de $f_1$ et $f_2$ sont triangulables.
L'application de Poincar\'e repr\'esente le changement de
coordonn\'ees locales rendant des it\'er\'es de $f_1$ et $f_2$
triangulaires. Cette application se prolonge holomorphiquement en
une application de $\C^k$ dans $\P^k$. L'image r\'eciproque de la
fonction de Green commune de $f_1$ et $f_2$ par l'application de
Poincar\'e est invariante par les applications
triangul\'ees.
\par    
Notons $[w_0:w_1:\cdots:w_k]$ les coordonn\'ees homog\`enes de
$\P^k$ et posons $z_s:=w_s/w_0$ pour $s=1,\ldots,k$.
Soit $f$ un endomorphisme holomorphe de degr\'e $d\geq
2$ de $\P^k$.
Un {\it relev\'e} de $f$ est une application polynomiale homog\`ene
$F:\C^{k+1}\longrightarrow \C^{k+1}$ v\'erifiant $F^{-1}(0)=\{0\}$ et
$\pi\circ f=F\circ\pi$ o\`u $\pi$ est l'application canonique de
$\C^{k+1}\setminus\{0\}$ dans $\P^k$. L'application $F$ est
d\'efinie \`a une constante multiplicative pr\`es. 
\par
La suite de
fonctions $d^{-n}\log\|F^n(w)\|$ converge vers une fonction
continue, plurisousharmonique $G$ qu'on appelle {\it la fonction de
Green}. Cette fonction v\'erifie $G\circ F=dG$ et $G(\lambda
z)=\log|\lambda|+G(w)$. De plus, toute fonction $v$ v\'erifiant
ces propri\'et\'es est telle que $v\leq G$. On peut d\'efinir un
courant positif ferm\'e $T$ de bidegr\'e $(1,1)$ de $\P^k$ par
la relation $\pi^*T:=dd^c G$.
C'est un courant de masse 1, invariant par $f$: $f^*
T=d.T$. 
\par
Pour tout $1\leq s\leq k$, 
on appelle {\it ensemble de Julia d'ordre $s$} le
support $J_s$ du courant $T^s:=T\wedge\ldots\wedge T$
qui est un courant positif ferm\'e de bidegr\'e $(s,s)$. Notons
$\mu:=T^k$ la mesure de probabilit\'e invariante de $f$, dite {\it
  mesure d'\'equilibre}.
Les ensembles de Julia ne sont pluripolaires dans aucun
ouvert qui les rencontre. 
\begin{theoreme}[Fornaess-Sibony \cite{Sibony}] 
Soit $f$ un endomorphisme holomorphe
  de $\P^k$ de degr\'e $d\geq 2$. La mesure $\mu$, dite mesure
  d'\'equilibre satisfaisant l'\'equation $f^*\mu=d^k\mu$. Il existe
  un ensemble pluripolaire $\E^*$ tel que pour $a\not\in \E^*$ les
  mesures 
$$\mu^a_n:=\frac{1}{d^{kn}}\sum_{f^n(a_i)=a}\delta_{a_i}$$
convergent vers $\mu$ o\`u on a not\'e $\delta_{a_i}$ la masse de Dirac en
$a_i$.
\end{theoreme}
\begin{theoreme}[Briend-Duval \cite{BriendDuval1}] Soit $f$ un
  endomorphisme holomorphe de $\P^k$ de degr\'e $d\geq 2$. Soit $A_n$
  l'ensemble des points p\'eriodiques r\'epulsifs d'ordre $n$ de
  $f$. Alors les mesures 
$$\frac{1}{d^{kn}}\sum_{a_i\in A_n}\delta_{a_i}$$
convergent vers $\mu$. En particulier, les points
p\'eriodiques r\'epulsifs de $f$ sont denses dans $J_k$.
\end{theoreme}
\par
Lorsque $f$ est polynomiale, c'est-\`a-dire si $f$ est
aussi un endomorphisme de $\C^k=\P^k\setminus\{w_0=0\}$,
on peut d\'efinir le
taux d'\'echappement vers l'infini des orbites de $f$
(not\'ee encore $G$ et appel\'e aussi {\it fonction de Green})
$G(z):=\lim_{n\rightarrow\infty} d^{-n}\log^+\|f^n(z)\|$. 
C'est
une fonction continue, plurisouharmonique, \`a croissance
logarithmique \`a l'infini (l'abus de
notation ne pr\^ete pas \`a confusion). Cette fonction est
\'egale \`a la restriction de la
fonction de Green d\'efinie pr\'ec\'edemment,
\`a l'hyperplan $\{w_0=1\}$ de $\C^{k+1}$ (ici l'application $F$
est choisie de sorte que sa premi\`ere fonction coordonn\'ee soit \'egale
\`a $w_0^d$). Dans ce cas, on a
$T=dd^c G$. 
Notons que les z\'eros de cette fonction sont les points d'orbite
born\'ee. En particulier, $G$ s'annule aux points p\'eriodiques.
\par  
Dans la suite de ce paragraphe,
on note $f_1$, $f_2$ deux endomorphismes
holomorphes permutables de degr\'es $d_1\geq 2$
et $d_2\geq 2$ de $\P^k$.
Soient $F_1$ et $F_2^*$ des relev\'es de $f_1$ et
$f_2$. La relation $f_1\circ f_2=f_2\circ f_1$ entra\^{\i}ne
$F_1\circ F_2^*=\lambda F_2^*\circ F_1$
avec $\lambda\not =0$.
Posons $F_2=\theta F_2^*$ avec $\theta^{d_1-1}=\lambda$.
L'application $F_2$ est un relev\'e de $f_2$. On
v\'erifie facilement que $F_1\circ F_2=F_2\circ F_1$. 
On note $G_1$, $G_2$ les fonctions de Green de $f_1$, $f_2$
d\'efinies gr\^ace \`a $F_1$ et $F_2$. 
\begin{proposition} Soient $f_1$ et $f_2$ deux endomorphismes
holomorphes permutables de $\P^k$. On a alors $G_1=G_2$.
En particulier, les ensembles
de Julia d'ordre $s$ de $f_1$ et $f_2$ sont \'egaux pour tout
$1\leq s\leq k$.
\end{proposition}
\begin{preuve} 
Consid\'erons les fonctions
$H_n:=d_2^{-n}G_1\circ F_2^n$. Ces fonctions v\'erifient
les propri\'et\'es suivantes:
\begin{enumerate}
\item $H_n\circ F_1=d_1H_n$.
\item $H_n(\lambda z)=\log|\lambda|+H_n(z)$ pour tout
$\lambda\not=0$.
\end{enumerate}
Par cons\'equent, $H_n\leq G_1$. La fonction
$L(z):=G_1(z)-\log\|z\|$ est
born\'ee car $G_1$ est continue et $G_1(\lambda z)=\log
|\lambda|+ G(z)$. On d\'eduit facilement de la relation
$H_n(z)=d_2^{-n}\log \|F_2^n(z)\|+ d_2^{-n}L\circ F_2^n(z)$ que $H_n$
tend vers $G_2$ lorsque $n$ tend vers l'infini. Par suite 
$G_2\leq G_1$. On montre de m\^eme que $G_1\leq G_2$.
\end{preuve}
\begin{lemme} Les endomorphismes $f_1$ et $f_2$ poss\`edent une
infinit\'e de points p\'eriodiques communs qui sont r\'epulsifs
pour $f_1$.
De plus, si $a$ est un point p\'eriodique r\'epulsif pour $f_1$,
alors il existe un entier $m\geq 0$ tel que $f_2^m(a)$ soit p\'eriodique
r\'epulsif pour $f_1$ et p\'eriodique pour $f_2$. 
\end{lemme}
\begin{preuve} Soit $a$ un point p\'eriodique r\'epulsif
d'ordre $n$ de $f_1$. Comme $f_1^n\circ f_2=f_2\circ
f_1^n$, $b:=f_2(a)$ est p\'eriodique d'ordre $s$ pour $f$ avec $s|n$. 
D'apr\`es le lemme 2.3 appliqu\'e \`a $f(z):=f_1^n(z+b)-b$,
$g(z):=f_1^n(z+a)-a$ et $h(z):=f_2(z+a)-b$, le point
$b$ est r\'epulsif
pour $f_1$. 
\par
Notons pour tout $n$, $A_n$ l'ensemble des
points p\'eriodiques r\'epulsifs d'ordre $n$ pour $f_1$.
Alors l'ensemble fini $A_1\cup A_2\cup\ldots\cup A_n$ est
invariant par $f_2$. Par cons\'equent, il existe un $m$ tel que
$f_2^m(a)$ soit p\'eriodique pour $f_2$. 
\par 
Pour tout nombre premier $n$, l'ensemble fini
$A_1\cup A_n$ est invariant par
$f_2$. Par cons\'equent, pour $n$ premier suffisamment grand,
$A_n$ est invariant par $f_2$. D'apr\`es le Th\'eor\`eme 3.2, 
$A_n$ est non vide pour $n$ suffisamment grand.
On en d\'eduit qu'il
existe au moins un point p\'eriodique
$a_n$ de $f_2$ appartenant \`a $A_n$. Il est clair que
l'ensemble de tels points est infini.
\end{preuve}
\par
Nous dirons qu'une application  
$\varphi:\C^k\longrightarrow \C^k$ est 
{\it rigide} s'il existe une application lin\'eaire
$l:\C^k\longrightarrow \C^{k-1}$ telle que pour tout $u\in\C^{k-1}$ 
l'image de $l^{-1}(u)$
par $\varphi$ soit une droite complexe passant
par $0$.
\begin{proposition} Soit $f$ un endomorphisme holomorphe
de degr\'e $d\geq 2$ de $\P^k$. Soit 
$b$ un point fixe r\'epulsif pour $f$. Alors il
existe une application holomorphe
$\varphi:\C^k\longrightarrow \P^k$
(dite l'application de Poincar\'e) et une application
triangulaire $\Lambda$ telles que
$\varphi(0)=b$, $\varphi'(0)$ inversible et 
$f\circ \varphi=\varphi\circ\Lambda$.
De plus, $\P^k\setminus \varphi(\C^k)$ est un ferm\'e
pluripolaire. Si $f$ est polynomiale on a
$G^*\circ \Lambda=dG^*$ o\`u
$G^*:=G\circ\varphi$. Si $f$ est polynomiale homog\`ene,
$\varphi$ est rigide. 
\end{proposition}
\begin{preuve} D'apr\`es le th\'eor\`eme de Sternberg, il existe
une application holomorphe $\varphi$ d'un voisinage de $0\in\C^k$
dans $\P^k$ telle que $\varphi(0)=b$, $\varphi'(0)$ 
inversible et telle que $\Lambda:=\varphi^{-1}\circ
f\circ\varphi$ soit triangulaire. On a \'egalement $\varphi\circ
\Lambda^n=f^n\circ\varphi$ pour tout $n\geq 1$. 
Comme $b$ est r\'epulsif pour
$f$, $\Lambda$ est dilatante. Par cons\'equent, l'application
$\varphi$ se prolonge en une application holomorphe de $\C^k$
dans $\P^k$ en posant $\varphi(z):=f^n\circ\varphi\circ\Lambda^{-n}(z)$
pour tout $z\in\C^k$ et pour un $n$ suffisamment grand. On
v\'erifie sans peine que la d\'efinition ne d\'epend pas de $n$.
\par
Il en r\'esulte que si $z\not\in\varphi(\C^k)$, l'ensemble
$\overline{\bigcup_{n\geq 0} f^{-n}(z)}$ ne contient pas $b$. En
particulier, la mesure $\mu^z_n$ du th\'eor\`eme 3.1
ne tend
pas vers $\mu$ quand $n\rightarrow+\infty$.
Par cons\'equent, $\P^k\setminus\varphi(\C^k)$ est pluripolaire; il
est ferm\'e car $\varphi$ est ouverte.
\par
Si $f$ est polynomiale on a $G\circ f=dG$. D'o\`u
$$G^*\circ\Lambda=G\circ\varphi\circ\Lambda=G\circ
f\circ\varphi=dG\circ\varphi=dG^*$$
\par
Supposons maintenant que $f$ soit polynomiale homog\`ene. Sans
perte de g\'en\'eralit\'e, on peut supposer que $b=(0,...,0,1)$.
Posons $f=(P_1,\ldots,P_k)$. L'application $h:=[P_1:\cdots:P_k]$
d\'efinit un endomorphisme holomorphe de $\P^{k-1}$ et
$b':=[0:\cdots:0:1]$ est un point fixe de $h$. Posons
$u_m:=z_m/z_k$ pour $m=1,\ldots,k-1$.
Alors $u':=(u_1,\ldots,u_{k-1})$
est un syst\`eme de coordonn\'ees affines de $\P^{k-1}$.
On a $h(u')=(R_1,\ldots,R_{k-1})$ au voisinage de
$b'$ o\`u $R_m:=P_m(u',1)/P_k(u',1)$. Posons $u_k:=d\log z_k
+\log P_k(u',1)/(d-1)$. On v\'erifie facilement que $u:=(u',u_k)$
est un syst\`eme de coordonn\'ees d'un voisinage de $b$. On a
\'egalement $f(u)=(h(u'),du_k)$. D'apr\`es la partie
pr\'ec\'edente appliqu\'ee \`a $h$ en $b'$, il existe une
application $\varphi_h:\ \C^{k-1}\longrightarrow \C^{k-1}$ telle
que $\varphi_h^{-1}\circ h\circ \varphi_h$ soit triangulaire.
Soit $\varphi:\C^k\longrightarrow \C^k$ d\'efinie par
$\varphi(z',z_k):=(\varphi_h(z'),z_k)$. Alors $\varphi^{-1}\circ
f\circ\varphi$ est triangulaire. Il est clair que $\varphi$ est
rigide pour l'aplication $l:\C^k\longrightarrow \C^{k-1}$ avec $l(z):=z'$. 
\end{preuve}
\begin{proposition} Soient $f_1$, $f_2$ deux endomorphismes
  holomorphes de $\P^k$ de fonction de Green commune $G$. 
Soit $a$ un point fixe commun \`a $f_1$,
$f_2$ r\'epulsif pour $f_1$. 
Soient $\lambda_1$, ..., $\lambda_k$ les valeurs propres de
$f_1'(a)$ rang\'ees dans l'ordre croissant de leurs modules. 
Alors il
existe une application holomorphe
$\varphi:\C^k\longrightarrow \P^k$
et des applications permutables $\Lambda_i\in \G_\lambda^*$ telles que
$\varphi(0)=a$, $\varphi'(0)$ inversible et 
$f_i\circ \varphi=\varphi\circ\Lambda_i$ o\`u
$\lambda:=(\lambda_1,\ldots,\lambda_k)$.
De plus, $\P^k\setminus \varphi(\C^k)$ est un ferm\'e
pluripolaire. Si les $f_i$ sont polynomiales
on a, $G^*\circ \Lambda_i=d_iG^*$ o\`u
$G^*:=G\circ\varphi$. Si les $f_i$ sont polynomiales homog\`enes,
$\varphi$ est rigide. 
\end{proposition}
\begin{preuve} D'apr\`es la proposition 2.2, il existe une
application holomorphe $\varphi$ d'un voisinage de $0\in\C^k$ \`a
valeurs 
dans $\P^k$ telle que $\varphi(0)=a$, $\varphi'(0)$ soit
inversible et $\varphi^{-1}\circ f_i \circ
\varphi$ soient $\lambda$-triangulaires. 
Posons
$\Lambda_i:=\varphi^{-1}\circ f_i\circ\varphi$. 
\par
Comme dans la proposition 3.5, on montre que $\varphi$ se prolonge
en une application holomorphe de $\C^k$ dans $\P^k$ et que
$\P^k\setminus \varphi(\C^k)$ est un ferm\'e pluripolaire. De plus
$G^*\circ\Lambda_i=d_iG^*$ si les $f_i$ sont polynomiales. Lorsque
les $f_i$ sont polynomiales homog\`enes, on peut choisir
comme dans la
proposition 3.5 une application $\varphi$ rigide.
\end{preuve}
\begin{remarque} \rm
Notons $T^*:=dd^cG^*$,
$\mu^*:=(T^*)^k$ et
$J_s^*$ le support de $(T^*)^s$ pour tout $1\leq
s\leq k$. On a
$J_s^*=\varphi^{-1}(J_s)$, $(T^*)^s
=\varphi^*(T^s)$, $\mu^*=\varphi^*(\mu)$,
$\Lambda_i^{\pm 1}(J_s^*)=J_s^*$
et $\Lambda_i^* ((T^*)^s)=d_i^s(T^*)^s$.
\end{remarque}
\section{Laminations de la mesure d'\'equilibre}
\ahead
Nous supposons dans ce paragraphe que $f_1$ et $f_2$ sont deux
applications polynomiales permutables de $\C^k$ qui se prolongent en
des endomorphismes holomorphes de $\P^k$.
Nous montrons sous l'hypoth\`ese du th\'eor\`eme 1.1 qu'il existe un ouvert
$\Omega\subset \C^k$ v\'erifiant les propri\'et\'es suivantes:
\begin{enumerate}
\item  $J_k\cap \Omega$ est une sous-vari\'et\'e r\'eelle analytique
non vide de $\Omega$.
\item Il existe une forme r\'eelle analytique $\phi$
de degr\'e maximal d\'efinie sur $J_k\cap \Omega$ telle que
$\mu=\phi\wedge [J_k]$ dans l'ouvert $\Omega$.
\end{enumerate}
\par
Sous l'hypoth\`ese de la proposition 3.6, on montrera que $a$ est
r\'eplusif pour $f_2$ et que le groupe ferm\'e engendr\'e par
les automorphismes $\Lambda_1$ et $\Lambda_2$ contient un
sous-groupe multiplicatif \`a un param\`etre r\'eel. Par
cons\'equent, l'ensemble $J^*:=\varphi^{-1}(J_k)$ est "lamin\'e" par
les courbes r\'eelles analytiques qui sont des orbites de ce
sous-groupe.
On obtient la premi\`ere partie gr\^ace \`a l'utilisation
des laminations obtenues par des points p\'eriodiques diff\'erents.
La deuxi\`eme partie
se d\'eduit \'egalement en utilisant la structure lamin\'ee de la mesure
$\mu$. Cette approche est celle d'Eremenko en dimension 1 \cite{Eremenko}.
\begin{lemme} Soit $\{g_n\}$ une famille non \'equicontinue
de germes d'applications
holomorphes d\'efinies au voisinage de $0\in\C^k$ \`a valeurs dans
$\C^m\subset\P^m$. Alors
il existe un vecteur non nul $v\in\C^k$ , une suite de nombres
complexes $\{z_i\}$ tendant vers $0$, une suite de nombres r\'eels
positifs $\{\rho_i\}$ tendant vers $0$ et une suite croissante
d'entiers
positifs $\{n_i\}$ tels que $g_{n_i}(z_iv+\rho_i \xi v)$ tende vers
une application holomorphe non constante $h(\xi)$ de $\C$ dans
$\C^m$. 
\end{lemme}
\begin{preuve} Montrons d'abord qu'il existe une droite $L$
passant par $0$ telle que la famille $\{g_{n|L}\}$ ne soit pas
\'equicontinue en $0$. Raisonnons par l'absurde. Supposons que pour
toute droite $L$ passant par $0$, la famille $\{g_{n|L}\}$ est
\'equicontinue en $0$. Notons ${\cal F}$ la famille des droites
complexes passant par $0$.
On peut supposer que
$g_n(0)=0$ pour tout $n$. Soit $U$ un voisinage suffisamment petit de
$0$. Alors pour toute droite
$L$ il existe $r_L$ rationnel positif tel que pour tout $n$
l'image de $L\cap\{\|z\|\leq r_L\}$ par $g_n$ soit contenue dans
$U$. Par cons\'equent, il existe un $r>0$ et
une famille non pluripolaire 
${\cal F}'\subset {\cal F}$ tels que pour tout $n$ et pour toute $L\in
{\cal F}'$ l'image de
$L_r:=L\cap\{\|z\|\leq r\}$ par $g_n$ soit contenue dans $U$.
Il existe une constante $M$ telle que $\|g_n(z)\|\leq M$ pour tout
$z\in L_r$ et pour tout $n$. D'apr\`es \cite{SibonyWong, Alexander},
l'enveloppe polynomialement convexe 
de $\overline{\bigcup_{L\in{\cal F}'}L\cap \{\|z\|\leq
r\}}$ contient un voisinage de $0$. Alors
il existe $r'>0$ tel que 
$\|g_n(z)\|\leq M$ pour tout $n$ et tout $|z|\leq
r'$. D'apr\`es le th\'eor\`eme de Montel, 
la famille $\{g_n\}$ est \'equicontinue
en $0$. 
C'est la contradiction
recherch\'ee.
\par
Soit $L$ une droite passant par $0$ telle que la famille
$\{g_{n|L}\}$ ne soit pas \'equicontinue en $0$. Soit $v\in L$ un
vecteur non nul. Posons $h_n(\xi):=g_n(\xi v)$. Alors la famille
$\{h_n\}$ n'est pas \'equicontinue en $0$. Il faut montrer
qu'il existe $\{n_i\}$, $\{z_i\}$, $\{\rho_i\}$ et $h$ tels que
$h_{n_i}(z_i+\rho_i \xi)$ tende vers $h$. Ceci dans le cas o\`u
$m=1$ est un lemme d\^u \`a Zalcman \cite{Schwick}. La preuve au cas
o\`u $m\geq 2$
se d\'eroule exactement de la m\^eme mani\`ere. 
\end{preuve}
\begin{lemme} Soit $A_c:=\{(m,n)\in \Z^2|\ 
d_1^md_2^n\leq c\}$ o\`u $c$ est une constante. 
Alors sous l'hypoth\`ese de la proposition 3.6,
la famille des $\Lambda_1^m\circ\Lambda_2^n$ avec $(m,n)\in A_c$
est un sous-ensemble born\'e de $\G_\lambda$.
\end{lemme}
\begin{preuve} On peut se limiter au cas o\`u $f_1$ et $f_2$ sont des
  endomorphismes polynomiaux. Supposons que les coefficients de
$\Lambda_1^m\circ\Lambda_2^n$ ne soient pas born\'es. Alors
l'ensemble de
telles applications n'est pas \'equicontinue. D'apr\`es le 
lemme 4.1, il existe des entiers relatifs $m_i$, $n_i$, un
vecteur non nul $v$, des suites $\{z_i\}$ et $\{\rho_i\}$ tendant
vers $0$ tels que $d_1^{m_i}d_2^{n_i}\leq c$ et tels que
$\Lambda_1^{m_i}\circ\Lambda_2^{n_i}(z_iv+\rho_i
\xi v)$ tende vers une application holomorphe non constante
$h(\xi)$ de
$\C$ dans $\C^k$. Comme $G^*$ est continue, on a
$$G^*(h(\xi))=\lim G^*(\Lambda_1^{m_i}\circ\Lambda_2^{n_i} (z_i
v+\rho_i \xi v))= \lim d_1^{m_i}d_2^{n_i}G^*(z_i v+\rho_i \xi v) =0.$$
En effet, $G^*(0)=0$, $d_1^{m_i}d_2^{n_i}\leq c$ et $z_i$, $\rho_i$
tendent vers $0$. La fonction $G$ s'annulle donc sur
$\varphi(h(\C))$. Or l'ensemble de
z\'eros de $G$ est un compact de $\C^k$; il ne peut donc
contenir une image holomorphe de $\C$.
\end{preuve}
\begin{corollaire} Soient $f_1$, $f_2$ deux endomorphismes
  satisfaisant l'hypoth\`ese de la proposition 3.6. 
Supposons de plus que les suites d'it\'er\'es
de $f_1$ et $f_2$ soient disjointes.
Alors le groupe
d'automorphismes polynomiaux de $\C^k$ engendr\'e par
$\Lambda_1$ et $\Lambda_2$ n'est pas discret. 
Le groupe ferm\'e engendr\'e par $\Lambda_1$
et $\Lambda_2$ contient un sous-groupe additif \`a un
param\`etre r\'eel $\{\Lambda^t\}_{t\in\R}$ v\'erifiant
$G^*\circ\Lambda^t=\exp(ct)G^*$ 
o\`u $c=1$ si $d_1^{n_1}\not=d_2^{n_2}$ pour tout 
$(n_1,n_2)\not=(0,0)$ et $c=0$ sinon. De plus si $c=1$ on a
$\lim_{t\rightarrow -\infty} \Lambda^t=0$.
\end{corollaire}
\begin{preuve} Soient $m_i$, $n_i$ des entiers relatifs
v\'erifiant $\lim d_1^{m_i}d_2^{n_i}=1$, $0<m_i<m_{i+1}$ et
$0>n_i>n_{i+1}$  pour
tout $i$. Les suites d'it\'er\'es de $f_1$ et $f_2$ \'etant
disjointes, les germes $f_1^{m_i}\circ f_2^{n_i}$, d\'efinis au
voisinage de $a$, sont
deux \`a deux diff\'erents. L'application de Poincar\'e \'etant la m\^eme,
les automorphismes $\Lambda_1^{m_i}\circ\Lambda_2^{n_i}$ sont
donc deux \`a deux diff\'erents. D'apr\`es le lemme 4.2, 
cette suite est born\'ee dans $\G_\lambda$, 
donc elle a des points d'accumulation et
n'est pas discr\`ete.
De m\^eme la suite $\Lambda_1^{-m_i}\circ\Lambda_2^{-n_i}$ est
auusi born\'ee. On en d\'eduit que  
les points d'accumulation de ces deux suites sont 
inversibles. Par cons\'equent, le groupe d'automorphismes engendr\'e
par $\Lambda_1$ et $\Lambda_2$ n'est pas discret.
\par
Notons $\Gamma$ le groupe ferm\'e engendr\'e par $\Lambda_1$ et
$\Lambda_2$. C'est un sous-groupe de $\G_\lambda^*$ qui
est un groupe de Lie commutatif de dimension finie. 
Par cons\'equent, $\Gamma$ est un sous-groupe de
Lie. Soit $\Gamma_0$ la composante de
l'identit\'e. On a $\dim \Gamma_0\geq 1$. Par cons\'equent, $\Gamma_0$
contient un sous-groupe \`a un param\`etre $(\Lambda^t)_{t\in \R}$.
\par
Consid\'erons l'application
$\Psi:\Gamma_0\longrightarrow \R$ d\'efinie
par la relation $G^*\circ\Lambda =
\exp(\Psi(\Lambda))G^*$ pour tout $\Lambda\in \Gamma_0$.
C'est un morphisme de
groupes de Lie. 
\par
Si le groupe $\Gamma_0$ est compact, l'ensemble 
$\Psi(\Gamma)$ est un sous-groupe discret de $\R$. Or l'image de $\Psi$
contient l'ensemble $\{m\log d_1+n\log d_2 \mbox{ avec } (m,n)\in \Z^2\}$.
Si cet ensemble est dense dans $\R$, c'est le cas si $\log d_2/\log
d_1$ est irrationnel, alors $\Psi$ est surjective et on peut choisir 
un sous-groupe $(\Lambda^t)_{t\in \R}$ sur lequel
$\Psi$ est surjective. 
\par
Posons
$\alpha(t):=\Psi(\Lambda^t)$. On a
$\alpha(t+t')=\alpha(t)\alpha(t')$. Il existe donc un $c\in
\R$ tel que $\alpha(t)=\exp(ct)$. Si
$d_1^{n_1}\not=d_2^{n_2}$ pour tout $(n_1,n_2)\not=(0,0)$, 
on a $\log d_2/\log d_1$
irrationnel et donc $c\not=0$. Le changement de
param\`etre $t\mapsto t/c$ permet de prendre $c=1$.
Dans le cas contraire, 
$\Psi$ n'est pas surjectif et n\'ecessairement $c=0$.    
\par
Supposons que $c=1$. Notons $\Lambda\in \G_\lambda$ une valeur
adh\'erente de la famille $\{\Lambda^t\}$ pour $t\rightarrow -\infty$.
Par continuit\'e de $G^*$, la relation $G^*\circ\Lambda^t=\exp(t)G^*$
implique que $G^*\circ\Lambda=0$. Comme $G^*$ ne peut s'annuler sur
aucune image holomorphe non constante de $\C$, on a n\'ecessairement
$\Lambda=0$. 
\end{preuve} 
\par
Nous allons introduire une notion de lamination adapt\'ee \`a nos
besoins.  
Soit $X$ une vari\'et\'e r\'eelle analytique de dimension $n$ et soit
$J\subset X$ un ferm\'e et
$b\in X$. On dit que $J$ est {\it $m$-lamin\'e} en $b$ s'il existe un
voisinage $U$ de $b$ muni des coordonn\'ees r\'eelles analytiques
$(x_1,\ldots, x_m,x_{m+1},\ldots, x_n)$
tels que $\Omega=U\times V$, $J\cap \Omega=U\times K$ o\`u $U$
(resp. $V$)
est un ouvert de $\R^m$ (resp. $\R^{n-m}$) et $K$ est un ferm\'e
de $V$.
Par d\'efinition, si $J$ est
$m$-lamin\'e en $b$, il est $m$-lamin\'e en tout point d'un voisinage
de $b$.
\par
Notons $\pi:\Omega=U\times V\longrightarrow V$ la projection de
$\Omega$ dans $V$. Soit $G$ une fonction
r\'eelle non n\'egative d\'efinie sur
$X$. On dit que dans $\Omega$, l'application
$\pi$ {\it lamine} $G$ s'il existe une fonction r\'eelle
positive $\delta$ d\'efinie
sur un intervalle $[0,\epsilon[$, tendant vers $1$ en $0$, telle que 
pour tout
$v\in V$ la restriction $G_v$ de $G$ dans
$U\times \{v\}$ soit r\'eelle
analytique et
$G_v(u)\leq \delta(\|u-u'\|)G_v(u')$ pour tous
$u, u'\in U$ et $\|u-u'\|<\epsilon$.
\par 
Soit $\mu$ une mesure de $X$. On dit que $\mu$ est {\it $m$-lamin\'ee}
en $b$ s'il existe un voisinage $\Omega=U\times V$ de $b$ muni des
coordonn\'ees r\'eelles analytiques locales
$(u,v)=(x_1,\ldots,x_m)\times (x_{m+1},\ldots, x_n)$ 
tel que $\mu$ soit le produit de la
mesure $dx_1\wedge\ldots\wedge dx_m$
d\'efinie sur $U$ et d'une mesure $\nu$
d\'efinie sur $V$. Par d\'efinition,
$\mu$ est $m$-lamin\'e en tout point de
$\Omega$. On dit que $G$ et $\mu$
sont {\it simultan\'ement $m$-lamin\'ees} si les
ouverts $\Omega$, $U$ et $V$ sont les
m\^emes pour $G$ et $\mu$. 
\par L'ensemble $J$ (resp. $\mu$) est {\it
$m$-lamin\'e(e)} dans un ouvert $W$ s'il
(resp. si elle) l'est en tout point
$b\in W$. 
\par
D'apr\`es le corollaire 4.3, les ensembles $J_s^*$ sont
$1$-lamin\'es en tout point sauf \'eventuellement en $0$. Par
cons\'equent, $J_s$ est $1$-lamin\'e en $b\in J_s\cap
\varphi(\C^k)$ lorsque $b\not=a$ et $\varphi^{-1}(b)$ n'est pas
contenu dans l'ensemble critique de $\varphi$. L'ensemble des
points $b$ qui ne v\'erifient pas ces deux propri\'et\'es, est un
ferm\'e pluripolaire ({\it voir} la proposition 3.5). 
\begin{proposition} Soient $f_1$ et $f_2$ deux endomorphismes
  polynomiaux permutables 
  de $\C^k$ qui se prolongent en des endomorphismes holomorphes de
  $\P^k$. Supposons que $d_1^{n_1}\not =d_2^{n_2}$ pour tous $n_1\geq
  1$, $n_2\geq 1$.
Alors il existe un ouvert
$\Omega=U\times V\subset \C^k$
muni des coordonn\'ees r\'eelles analytiques  $(x_1,\ldots,x_m)\times
(x_{m+1},\ldots, x_{2k})$ tel que $J_k\cap \Omega=U\times \{0\}$,
$\mu=dx_1\wedge\ldots\wedge dx_m\wedge [J_k]$
dans $\Omega$. De plus la
projection $\pi$ de $\Omega$ dans $V$ lamine la fonction $G$ et l'entier
  $m$ v\'erifie $k\leq m<2k$.  
\end{proposition}
\begin{preuve} 
Soit $m$ l'entier maximal tel qu'il existe 
un ouvert $\Omega=U\times V$, muni des
coordonn\'ees r\'eelles analytiques
$(x_1,\ldots,x_m)\times (x_{m+1},\ldots, x_{2k})$,
une mesure $\nu$ de
$V$ tels que $J_k\cap \Omega\not=\emptyset$ avec
$\mu=dx_1\wedge \ldots\wedge dx_m \otimes\nu$
dans $\Omega$, la
projection de $\Omega$ dans $V$ lamine $G$. 
Notons $K$ le support de $\nu$.
On a $J_k\cap \Omega=U\times K$.
D'apr\`es le th\'eor\`eme 3.2, il existe   un point
$b\in J_k\cap W$ p\'eriodique r\'epulsif
pour $f_1$. 
Pour simplifier les notations, supposons que $b=(0,0)$. Si
$0$ est un point isol\'e dans $K$, quitte \`a remplacer $\Omega$
par un ouvert convenable on peut supposer que $K=\{0\}$; 
la proposition est alors vraie. En effet, $m\geq k$ car la vari\'et\'e
r\'eelle analytique $J_k\cap
\Omega$ n'est pas pluripolaire et $m<2k$ car $f_1$ et $f_2$ sont
polynomiales et le support de $\mu$ ne peut contenir un ouvert
\cite[p.163]{Sibony}. 
\par
Supposons maintenant
que $0$ ne soit pas isol\'e dans $K$. Nous en d\'eduirons que $m$ n'est
pas maximal. L'id\'ee est qu'on peut fabriquer une $(m+1)$-lamination
\`a partir d'une $m$-lamination et d'une 1-lamination transverse \`a la
premi\`ere. La 1-lamination sera choisie gr\^ace au corollaire 4.3.
\par
D'apr\`es le lemme 3.4, 
pour un $n$ convenable le point
$a:=f_2^n(b)$ est un point p\'eriodique commun pour $f_1$ et
$f_2$, r\'epulsif pour $f_1$. 
Pour simplifier les notations, on
suppose que $a$ est fixe pour $f_1$ et $f_2$.
D'apr\`es le corollaire 4.3, on a $(\Lambda^t)^*(
\mu^*)=\exp(kt)\mu^*$ et $G^*\circ\Lambda^t=
\exp(t)G^*$. Soient $\Sigma$
l'ensemble critique de $f_2^n$, $Y:=\varphi^{-1}\circ f_2^n(\Sigma)$ et
$\Theta$ un voisinage suffisamment petit de $0$. On a 
$J^*\cap \Theta\not\subset Y$ car $Y$ est pluripolaire o\`u
$J^*:=\varphi^{-1}(J_k)$. Soient
$c\in \Theta\setminus Y$ et $z\in \Omega$ v\'erifiant
$f_2^n(z)=\varphi(c)$. Alors la mesure $\mu$ et la fonction $G$ sont
simultan\'ement $m$-lamin\'ees en $\varphi(c)$ 
car $f_2^n$ r\'ealise une biholomorphisme d'un voisinage de
$z$ dans un voisinage de $\varphi(c)$ satisfaisant les relations
$(f_2^n)^*(\mu)=d_2^n\mu$ et $G\circ f_2^n=d_2^n\circ G$.
Par cons\'equent,
$\mu^*$ et $G^*$ sont simultan\'ement $m$-lamin\'ees en $c$.
De plus, lorsque $c$
appartient \`a $J^*$, la mesure $\mu^*$ et la fonction $G^*$ 
ne sont pas simultan\'ement
$(m+1)$-lamin\'ees en $c$ car on a choisi $m$ maximal.
\par
Soit $x''\in K\setminus\{0\}$ un point suffisamment proche de $0$
v\'erifiant $U\times\{x''\}\not\subset \Sigma$. Posons
$U_0:=\varphi^{-1}\circ f_2^n(U\times\{0\})\cap \Theta$ et
$U_1:=\varphi^{-1}\circ f_2^n(U\times\{x''\})\cap \Theta$.
D'apr\`es le corollaire 5.14, l'orbite de tout point $c\in U_1$ par
$\Lambda^t$ est
une courbe r\'eelle analytique de limite $0$ quand $t\mapsto 0$. Par
cons\'equent, l'orbite d'un point g\'en\'erique $c\in U_1$ coupe
$U_1$ transversalement en $c$. On fixe un point
$c\in J^*$ v\'erifiant cette
propri\'et\'e. On montrera que $\mu^*$ et $G^*$ sont simultan\'ement
$(m+1)$-lamin\'ee en $c$, ce sera la contradiction cherch\'ee. 
\par
On choisit
un syst\`eme de coordonn\'ees r\'eelles analytiques
$(y_1,\ldots, y_m,y_{m+1},\ldots,y_{2k})$
d'un voisinage $\Omega_2=U_2\times V_2$ de $c$ et une mesure 
$\nu^*$ de $V_2$ tels que $c=(0,0)$ et
$\mu^*=dy_1\wedge \ldots
\wedge d y_m\otimes \nu^*$ dans $\Omega_2$ et
tels que la projection $\pi_1$ de
$\Omega_2$ dans $V_2$ lamine $G^*$.
Pour tout $z\in \Omega_2$, notons
$\gamma_z$ l'orbite de $z$ par $\Lambda^t$.
On va montrer qu'apr\`es changement de coordonn\'ees, on peut supposer
que $\{0\}\times V_2$ est lamin\'e par les courbes $\gamma_z$; un
second changement de coordonn\'ees les redresse et fournit une
$(m+1)$-lamination.
\par
Sans perte de g\'en\'eralit\'e,
on peut supposer que la droite tangente de $\gamma_c$ en $c$ est
$\{y_1=\cdots=y_m=y_{m+2}=\cdots=y_{2k}=0\}$. Posons
$L:=\{y_1=\cdots=y_{m+1}=0\}\cap \Omega_2$ et $H:=\bigcup_{z\in
L}\gamma_z\cap \Omega_2$. Pour $\Omega_2$ suffisamment petit, $H$
est le graphe d'une application r\'eelle analytique
$h=(h_1,\ldots,h_m)$
au-dessus de $V_2$.  Quitte \`a effectuer le changement
de coordonn\'ees
$(y_1,\ldots,y_{2k})\mapsto (y_1-h_1,\ldots,y_m-h_m,y_{m+1},
\ldots,y_{2k})$ on
peut supposer que $h=0$. On identifie $H$ \`a $\{0\}\times V_2$.
Il suffit maintenant de
montrer que $\nu^*$ et $G^*_{|H}$ sont
simultan\'ement 1-lamin\'ees. D'apr\`es le
corollaire 4.3, $G^*_{|H}\circ\Lambda^t=\exp(t) G^*_{|H}$ et 
$(\Lambda^t)^*(\nu_H)=\exp(kt)\nu_H$ dans
$H$ pour $|t|$ assez petit. On peut
trouver un syst\`eme de coordonn\'ees r\'eelles analytiques
$(u,y_{m+2}',\ldots, y_{2k}')$ de $H$ tel que
$G^*_{|H}\circ\tau=\exp(t)G^*_{|H}$ et
$\tau_t^*(\nu^*)= \exp(kt)\nu^*$ dans
$H$ pour $|t|$ assez petit o\`u
$\tau_t(u,y_{m+2}',\ldots,y_{2k}'):=(u+t,y_{m+2}',
\ldots, y_{2k}')$. Posons
${\nu^*}':=\exp(-ku)\nu^*$.
Alors $\tau_t^*({\nu^*}')={\nu^*}'$ dans $V_2$.
Par cons\'equent, il existe une mesure
$\sigma$ de $\{u=0\}$ telle que
${\nu^*}'=du\otimes \sigma$.
D'o\`u $\nu^*=\exp(ku)du\otimes \sigma$. Posons
$y_{m+1}':=k^{-1}\exp(ku)$.
Dans les coordonn\'ees $(y_{m+1}',\ldots,y_{2k}')$
de $H$, on a $\nu^*=dy_{m+1}'\otimes \sigma$.
On v\'erifie facilement
que la projection $( y_{m+1}',\ldots,y_{2k}')\mapsto
(y_{m+2}',\ldots,y_{2k}')$ lamine $G^*_{|H}$.
On en d\'eduit que $G^*$ et
$\mu^*$ sont simultan\'ement $(m+1)$-lamin\'ees
dans $\Omega_2$ par l'application 
$$(y_1,\ldots,y_m,y_{m+1}',\ldots,y_{2k}')\mapsto
 (y_1,\ldots,y_m,y_{m+1}').$$
Ceci contredit
le fait que $m$ est maximal. 
\end{preuve} 
\section{Lin\'earisation et groupe $\A$}
\ahead
Dans ce paragraphe on se propose de d\'emontrer le th\'eor\`eme
suivant qui dans le cas o\`u $m=2k-1$ et $f$ est homog\`ene, se
r\'eduit au th\'eor\`eme de Berteloot-Loeb \cite{BertelootLoeb}:
\begin{theoreme} Soit $f$ un endomorphisme polynomial de degr\'e
  alg\'ebrique $d\geq 2$ de $\C^k$
qui se prolonge en un endomorphisme holomorphe de $\P^k$. Soient
$G$, $\mu$ et $J_k$ sa fonction de Green, sa mesure d'\'equilibre et
son ensemble de Julia d'ordre maximal. Supposons qu'il existe
un ouvert $\Omega=U\times V\subset \C^k$
muni de coordonn\'ees r\'eelles analytiques $(x_1,\ldots,x_m)\times
(x_{m+1},\ldots, x_{2k})$ v\'erifiant $J_k\cap \Omega=U\times \{0\}$,
$\mu=dx_1\wedge\ldots\wedge dx_m\wedge [J_k]$
dans $\Omega$ et tels que la
projection $\pi$ de $\Omega$ dans $V$ lamine la fonction $G$.  
Alors il
existe une application holomorphe $\varphi: \C^k\longrightarrow
\C^k$ et une application affine holomorphe
$\Lambda$ satisfaisant l'\'equation $f\circ\varphi=\varphi\circ \Lambda$
et telles que $\varphi(\C^k)$
soit un ouvert de compl\'ement pluripolaire de $\P^k$. De plus, il
existe un groupe discret d'applications affines holomorphes
$\A$ de $\C^k$ agissant transitivement
sur les fibres de $\varphi$.
\end{theoreme}
\begin{remarque} \rm
En fait, le th\'eor\`eme est vrai en supposant que
  les coordonn\'ees de la lamination sont de classe ${\cal
  C}^2$. Quant \`a l'ensemble $J_k$, il est r\'eel analytique
  lorsqu'il est de classe ${\cal C}^2$. En effet, il est invariant par
  des applications localement triangulables.
Nous nous contentons d'en donner la preuve dans le cas r\'eel
  analytique car c'est ce cas que nous utilisons.
On verra que gr\^ace \`a la proposition 4.4, le th\'eor\`eme 1.1 sera
prouv\'e de la m\^eme mani\`ere. 
\end{remarque}
\par
Esquissons les id\'ees de la preuve.
Soit $b\in J_k\cap U$ un point p\'eriodique r\'epulsif de $f$. Dans un
premier temps, on suppose que $b$ est fixe pour simplifier les
notations. Le fait que la p\'eriode de $b$ ne soit
pas n\'ecessairement \'egale \`a 1 sera \'etudi\'e 
\`a la fin de la d\'emonstration du
th\'eor\`eme. 
\par
Notons
$\varphi:\C^k\longrightarrow \C^k$  une application de Poincar\'e de
$f$ en $b$ v\'erifiant
$\varphi(0)=b$, $\varphi'(0)$ inversible et telle que
$\Lambda:=\varphi^{-1}\circ f\circ \varphi$ soit triangulaire 
({\it voir} la proposition 3.5). On peut
supposer que les \'el\'ements de la diagonale principale de $\Lambda$
sont rang\'es selon les modules croissants.
On montrera que $\Lambda$ est lin\'eaire et d\'efinie par une
matrice diagonale. De plus, toute valeur propre de $\Lambda'(0)$
est \'egale \`a $\pm d$ ou \'egale \`a $\sqrt{d}$ en module.
La nature des valeurs propres de $\Lambda$ permettra, en utilisant
l'invariance de $J_k$ par $f$ de montrer que $J^*:=\varphi^{-1}(J_k)$
admet pour \'equations 
$$\Im z''=q(z',\overline z')$$
dans des coordonn\'ees convenables $z':=(z_1,\ldots,z_n)$,
$z'':=(z_{n+1},\ldots,z_k)$ et $z:=(z',z'')$, o\`u $q$ est une
application 
polynomiale homog\`ene de degr\'e 2. En changeant de coordonn\'ees on
peut \'eliminer les termes harmoniques de $q$. On montre ensuite que
$G^*:=G\circ\varphi$ est constante sur les vari\'et\'es
$$\Im z''=q(z',\overline z')+c, \  c\in \R^{k-n}.$$
Cela permet de prouver que les
biholomorphismes locaux qui permettent de passer d'un point sur la
fibre de $\varphi$ \`a un autre se prolongent en applications affines.   
Ces applications affines forment le groupe $\A$.
On v\'erifie ensuite que $\A$ op\`ere
transitivement sur les fibres de $\varphi$. Pour cela, on utilise 
la dynamique de $f$ sur $J_k$.
\par
Commen\c cons par quelques remarques sur la g\'eom\'etrie de $J^*$.
On sait que l'ensemble $J_k\cap \Omega$
n'est pas pluripolaire, or c'est une
vari\'et\'e r\'eelle analytique,
sa dimension $m$ est donc sup\'erieure ou \'egale \`a $k$. Posons $n:=m-k$. 
Le sous-espace tangent complexe de $J_k$ en
un point g\'en\'erique $z\in J_k\cap \Omega$ est de dimension
$n$ car sinon $J_k$ serait contenu dans une hypersurface complexe de
$\Omega$. Par continuit\'e, on peut choisir
$\Omega$ de sorte que $J_k\cap
\Omega$ soit {\it CR-g\'en\'erique}, c'est-\`a-dire le sous-espace
tangent complexe en tout
point de $J_k\cap \Omega$
soit de dimension $n$. Cela \'equivaut \`a dire que l'espace complexe
engendr\'e par l'espace tangent de $J_k\cap \Omega$ en tout point 
est de dimension $k$.
Comme $f$ est polynomiale, $J_k$ qui est
diff\'erent de $\C^k$, ne
contient aucun ouvert de $\C^k$ \cite[p.163]{Sibony} et par suite
$n<k$. 
\par
Notons
$H$ l'espace tangent r\'eel \`a $J^*$ en $0$ et 
$L$ l'espace tangent complexe de $J^*$ en
$0$. Dans un voisinage de $0$, $J^*$ est une sous-vari\'et\'e
r\'eelle analytique. Comme $J^*$ est invariant par l'automorphisme
dilatant $\Lambda$, c'est une sous-vari\'et\'e r\'eelle
analytique de $\C^k$.
L'ensemble $J^*$ \'etant invariant par $\Lambda$, il en r\'esulte que 
$H$ et $L$ sont invariants par $\Lambda'(0)$. Quitte \`a
effectuer un changement lin\'eaire de coordonn\'ees, on peut
supposer qu'il existe $1\leq s_1 \leq \cdots \leq s_n\leq k$ tels que
$H=\{\Im z_s=0 \mbox{ pour tout } s\not=s_1,\ldots,s_n\}$. 
Ce changement ne modifie pas la forme triangulaire de $\Lambda$.
Il est clair aussi que si $f$ est homog\`ene,
ce changement de coordonn\'ees ne change pas la rigidit\'e de
$\varphi$. 
\begin{lemme} Soient $\gamma_1$, ..., $\gamma_k$ les \'el\'ements de
  la diagonale principale de $\Lambda'(0)$. On a 
$|\gamma_s|\leq d$ pour tout $1\leq s\leq k$.
\end{lemme}
\begin{preuve} On a suppos\'e ci-dessus que
  $1<|\gamma_1|\leq\cdots\leq |\gamma_k|$. Il suffit de montrer que
  $|\gamma_k| \leq d$. 
Notons $l$ la droite r\'eelle $\{z_1=\cdots=z_{k-1}=\Im z_k=0\}$.
On a $l\subset H$ et la droite complexe engendr\'ee par $l$ 
est invariante par
  $\Lambda$ car $\Lambda$ est triangulaire. 
Soit $l'$ une courbe r\'eelle analytique contenue dans $J^*$ et
tangente \`a $l$. Soit $\{n_i\}$ une suite croissante d'entiers
  positifs telle que $\lim \arg \gamma_k^{n_i}\longrightarrow 0$.
Utilisant le d\'eveloppement de
Taylor des \'equations d\'efinissant $l'$, on v\'erifie
facilement que $\Lambda^{n_i}(l')$ tend vers $l$ quand $i$ tend vers
l'infini. La droite $l$ est
donc contenue dans $J^*$ car  $J^*$ est un ferm\'e invariant
par $\Lambda$. Notons $G_k$ la restriction de $G^*$ \`a la
droite complexe
$\{z_1=\cdots=z_{k-1}=0\}$. La fonction $G_k$ est 
sous-harmonique non identiquement nulle et elle ne prend que des
valeurs positives ou nulles. De plus, $G_k$ s'annule sur $\Im
z_k=0$ et $G_k(\gamma_k z)=dG_k(z)$. 
Donc on a, si $d<|\gamma_k|$, 
$$G_k(z)\leq \left(\frac{d}{|\gamma_k|}|z|+\const\right).$$
En effet on peut supposer que $G$ est radiale, 
$|\gamma_k|^s\leq |z|<|\gamma_k|^{s+1}$ et on obtient 
\begin{eqnarray*}
G(z) & \leq &
G(|\gamma_k|^{s+1})=d^{s+1}G(1)=G(1)\left(\frac{d}{|\gamma_k|}\right)^s
|\gamma_k|^s d \\
& \leq & |\gamma_k|^s d+\const\leq
|z|\frac{d}{|\gamma_k|}+\const
\end{eqnarray*}
Le principe de
Phragm\'en-Lindel\"of 
impliquerait que $G_k$ serait identiquement nulle. Par
cons\'equent, $|\gamma_k|\leq d$.   
\end{preuve}
\begin{proposition} Soit $H=\{\Im z_s=0 \mbox{ pour }
  s\not=s_1,\ldots,s_n\}$ l'espace tangent de $J^*$ en $0$. Alors
 $s_j=j$ pour tout $j=1,\ldots,n$. De
plus, $|\gamma_j|=\sqrt{d}$ pour $1\leq j\leq n$ et  
$\gamma_j=\pm d$ pour $n+1\leq j \leq k$. En particulier, $\Lambda$
est de degr\'e au plus deux et les \'equations de $J^*$ au voisinage
de $0$ sont de la forme 
$$\Im z_j=h_j(z_1,\overline z_1,\ldots, z_n,
\overline z_n), \ n+1\leq j\leq k$$
o\`u les $h_j$ sont des polyn\^omes homog\`enes de degr\'e 2.
\end{proposition}
\begin{preuve} Soit $s$ le plus grand indice v\'erifiant
$|\gamma_s|\leq \sqrt{d}$. Posons $s=0$ si
$|\gamma_1|>\sqrt{d}$. Posons 
$H_s:=\{z_1=\cdots=z_s=0\}$. 
Comme dans le lemme 5.3,
en utilisant le d\'eveloppement de Taylor des \'equations 
d'une vari\'et\'e contenue dans $J^*$ et tangente \`a
$H_s$ et en prenant des limites des images par $\Lambda^n$, 
on peut montrer que $J^*$ contient
$H\cap H_s$. Donc $H\cap H_s$ ne contient aucune
droite complexe. On en d\'eduit que $s_j\leq s$
pour tout $1\leq j\leq n$.
Donc $|\gamma_{s_j}|\leq \sqrt{d}$ pour tout $1\leq j\leq n$. Par
cons\'equent, le jacobien r\'eel de la d\'eriv\'ee de
$\Lambda_{|J^*}$ en $0$ est major\'e par $(\sqrt{d})^{2n}
d^{k-n}=d^k$. De la relation $f^*\mu=d^k\mu$ on d\'eduit que 
$\Lambda^*(\mu^*)=d^k\mu^*$ 
au voisinage de $0$ car $\varphi_*(\mu^*)=\mu$.
Le jacobien de $\Lambda_{|J^*}$ en $0$ est donc \'egal \`a
$d^k$. Par suite, $|\gamma_{s_j}|=\sqrt{d}$ pour tout
$1\leq j\leq n$ et $|\gamma_j|=d$ si $j\not=s_1,\ldots, s_n$.
En particulier, $\{1,\ldots,n\}=\{s_1,\ldots,s_n\}$. D'o\`u
$s_j=j$ pour $j=1,\ldots,n$.
\par
On a alors $H=\{\Im z_{n+1}=\cdots =\Im z_k=0\}$. Comme $H$ est
invariant par $\Lambda'(0)$, $\gamma_j$ est un nombre
r\'eel pour tout $j\geq n+1$. D'o\`u $\gamma_j=\pm d$ pour
$j\geq n+1$. 
\par
Le fait que $\Lambda$ est de degr\'e au plus deux
r\'esulte des r\'esonnances possibles entre les valeurs propres de
$\Lambda'(0)$. 
\par
La vari\'et\'e $J^*$ \'etant r\'eelle analytique et tangente \`a $H$ en
$0$, au voisinage de $0$ elle est d\'efinie par les
\'equations 
$$\Im z_j=h_j(z_1,\overline z_1,\ldots, z_n,
\overline z_n,\Re z_{n+1},\ldots, \Re z_k)$$
les fonctions $h_j$ \'etant r\'eelles analytiques
pour $j=n+1,\ldots,k$.  Utilisant le d\'eveloppement de Taylor des
$h_j$ et l'invariance de $J^*$ par $\Lambda$, on montre que les $h_j$
sont de la forme
$$h_j(z_1,\overline z_1,\ldots, z_n,
\overline z_n,\Re z_{n+1},\ldots, \Re z_k)= 
P_j(z_1,\overline z_1,\ldots, z_n,\overline z_n)+L_j(\Re
z_{n+1},\ldots, \Re z_k)$$
o\`u les $P_j$ sont des polyn\^omes homog\`enes de degr\'e 2 et les
$L_j$ sont lin\'eaires. La vari\'et\'e $J^*$ \'etant tangente \`a $H$,
les $L_j$ sont donc nuls.
\end{preuve}
\par
Nous avons besoins du lemme suivant que nous appliquerons \`a des
restrictions convenables de $G^*$. Pour tout $R>0$ on pose 
$$v_R(z):= 2R-\frac{2R}{\pi}\left(\arctan\frac{R-\Re z}{\Im z}+
\arctan\frac{R+\Re z}{\Im z}\right).$$
La somme dans les parenth\`eses repr\'esente l'angle en $z$ du
triangle de sommets $z$, $-R$ et $R$.
\begin{lemme} Soit $v\geq 0$ une fonction continue
sous-harmonique d\'efinie sur le demi-plan
$H^+:=\{z\in \C|\ \Im z\geq 0\}$ et
nulle sur $\R$. Soit $R$ un
nombre r\'eel positif.
Supposons qu'il existe
une constante $c\geq 0$ telle que $v(z)\leq c|z|$ pour tout
$z$ v\'erifiant $|z|=R$.
Alors $v(z)\leq c v_R(z)$ pour tout $z\in H^+$ v\'erifiant $|z|<R$.
S'il existe un nombre r\'eel $d>1$ tel que $v(dz)=dv(z)$ pour
tout $z\in H^+$, alors $v(z)=c\Im z$ o\`u $c\geq 0$ est
une constante.
\end{lemme}
\begin{preuve} 
La fonction $cv_R(z)$ est harmonique dans le demi-disque $\{|z|<R\}\cap
H^+$ nulle sur $]-R,R[\cap H^+$ et \'egale \`a
$cR$ sur $\{|z|=R\}\cap H^+$. Par cons\'equent, 
$cv_R$ majore la fonction $v$ dans ce demi-disque. D'o\`u $v(z)\leq
c v_R(z)$ pour tout $z\in H^+$ v\'erifiant $|z|<R$.
\par
Supposons que $v(dz)=dv(z)$. Soit $a\in H^+$ tel que $|a|=1$ et
$v(a)=\max_{\{|z|=1|\}\cap H^+}v(z)$. Posons $c=v(a)$.
On a $v(d^s z)\leq cd^s$ pour tout $s\geq 1$ et tout $|z|=1$. 
On d\'eduit de la
partie pr\'ec\'edente que $v(z)\leq c v_{R_s}(z)$ pour tout $z\in H^+$ et
tout $R_s:=d^s>|z|$.
Alors $v(z)\leq\lim_{s\rightarrow
\infty} c v_{R_s}(z)=c\Im z$. 
Pour $z=a$, cette in\'egalit\'e nous donne $a=i$. 
La fonction sous-harmonique, $v(z)-c\Im z$,  atteint donc son maximum en
$a$; elle est donc identiquement nulle et $v(z)=c\Im z$.
\end{preuve}
\begin{proposition} L'automorphisme $\Lambda$ est lin\'eaire et
d\'efini par une matrice diagonale. On peut choisir les
coordonn\'ees $z$ telles que  $J^*$ soit d\'efini par des
\'equations de la forme $\Im z''=q(z',\overline z')$ o\`u
$z':=(z_1,\ldots,z_n)$, $z'':=(z_{n+1},\ldots, z_k)$ et $q$ est une
forme hermitienne \`a valeurs vectorielles
v\'erifiant $q^{-1}(0)=\{0\}$. De plus les nouvelles coordonn\'ees ne
changent pas la rigidit\'e \'eventuelle de $\varphi$.
\end{proposition}
\begin{preuve} L'id\'ee est que l'existence de termes non diagonaux
  dans $\Lambda$ et
  l'invariance de $G^*$ par $\Lambda$ permettent de construire une
  droite complexe sur laquelle $G^*$ est nulle; ce qui est impossible.
\par
Notons $\Lambda_1$, ..., $\Lambda_k$ les fonctions
coordonn\'ees de $\Lambda$. Montrons d'abord que si $j\geq n+1$,
$\Lambda_j$ est ind\'ependante de $z_s$ pour tout $s\geq n+1$ et
$s\not=j$. 
On sait que 
$H_n:=\{z_1=\cdots=z_n=0\}$ est invariant par la d\'eriv\'ee 
$\Lambda'(0)$ et que
$\Lambda'(0)$ est
\'egale \`a $\Lambda$ sur $H_n$ car $\Lambda$ est triangulaire et
qu'il n'y a pas de r\'esonnances dans cet espace. Par
cons\'equent, $H_n$ est invariant par $\Lambda$.
D'apr\`es la proposition 5.4, 
$J^*\cap H_n=\{\Im z_{n+1}=\cdots=\Im
z_k=0\}$. 
L'application lin\'eaire $\Lambda_{|H_n}$ pr\'eserve le
sous-espace r\'eel $J^*\cap H_n$. Par
cons\'equent, ses coefficients sont r\'eels. Montrons 
que $\Lambda_{|H_n}$ est d\'efinie par une
matrice diagonale. Si tel n'\'etait pas le cas,
quitte \`a effectuer un changement lin\'eaire de
coordonn\'ees, on peut suposer que 
$\Lambda_{|H_n}$ contient le bloc de Jordan
$(\Lambda_{k-1},\Lambda_k)=(\alpha 
z_{k-1},\alpha z_k+d z_{k-1})$ o\`u $\alpha=\pm d$.
Quitte \`a remplacer
$\Lambda$ par $\Lambda^2$, on peut
supposer que $\alpha=d$. Posons
$K:=\{z_1=\cdots=z_{k-2}=0\}$.  La fonction $G^*$ \'etant continue, il
existe donc une constante $c>0$ telle que
$G^*(z)\leq c$ pour tout $z$ v\'erifiant $\|z\|\leq 1$. Soit $G_K$
la restriction de $G^*$ sur $K$.
Posons $\D:=\{z_{k-1}=0\}\cap K$,
$\D':=\{z_k=0\}\cap K$ et $\D_s:=\Lambda^s(\D')=\{(a,sa)\in K \mbox{
  avec } a\in \C\}$ pour tout $s\geq 1$. 
Les relations $G^*\circ \Lambda=d G^*$ et
$G^*(z)\leq c$ pour $\|z\|\leq 1$ impliquent $G_K(z)\leq c|z_k|/s$
pour tout $z=(z_{k-1},z_k)\in \D_s$ v\'erifiant $|z_{k-1}|\leq
d^s$. Puisque $G^*$ est nulle sur l'ensemble $J^*$ qui contient 
$\{\Im z_{k-1}=\Im z_k=0\}\cap K$,
d'apr\`es le lemme 5.5, pour tout $z=(z_{k-1},z_k)\in \D_s$ et
tout $s$ tel que $d^s\geq |z_{k-1}|$ on a 
$$G_K(z)\leq cv_{d^s}(z_{k-1})= cv_{d^s}(z_k/s).$$ 
Par continuit\'e on a
$$G_K(0,z_k)=\lim_{s\rightarrow \infty}G_K(z_k/s,z_k)\leq
\lim_{s\rightarrow \infty} cv_{d^s}(z_k/s)=0.$$
La fonction $G^*$ est donc nulle sur la droite
$\D$. C'est la contradiction cherch\'ee. Par suite
$\Lambda_{|H_n}$ est lin\'eaire et d\'efinie par une matrice
diagonale.
\par
Montrons maintenant que $\Lambda_j$ est ind\'ependant de $z_s$
pour tout $s\not=j$ et tout $1\leq j\leq n$. Raisonnons par l'absurde.
Supposons qu'il existe $1\leq j\leq n$
et $s\not=j$ tel que $\Lambda_j$
d\'epende de $z_s$. D'apr\`es la proposition 5.4, on a $s\leq j-1$ et 
$\Lambda_i$ est lin\'eaire pour tout $i\leq n$.
Quitte \`a faire un changement
lin\'eaire des coordonn\'ees en $z_1$, ..., $z_n$,
on peut supposer que
$\Lambda_n(z):=\gamma_n z_n +z_{n-1}$ et
$\gamma_{n-1}=\gamma_n$. Ce changement ne modifie pas la forme
triangulaire de $\Lambda$. 
Notons $H_{n-2}:=\{z_1=\cdots=z_{n-2}=0\}$.
Ce sous-espace est invariant par
$\Lambda$ car $\Lambda$ est triangulaire. 
L'ensemble $J^*\cap H_{n-2}$ est une
vari\'et\'e r\'eelle analytique dont
le plan tangent en $0$ est  $H\cap H_{n-2}$.
Cette vari\'et\'e est d\'efinie
par $k-n$ \'equations $\Im z_j=l_j(z_{n-1},
\overline z_{n-1},z_n,\overline z_n)$ o\`u 
$$l_j(z_{n-1},
\overline z_{n-1},z_n,\overline z_n):=h_j(0,\ldots,0,z_{n-1},
\overline z_{n-1},z_n,\overline z_n)$$
pour $j=n+1, \ldots,k$. 
%
%
\par
Fixons un $j\geq n+1$.
%
%
Quitte \`a
remplacer $\Lambda$ par $\Lambda^2$, on peut
supposer que $\gamma_j=d$. On a
${\Lambda_j}_{|H_{n-2}}(z)= d z_j+P(z_{n-1},z_n)$
o\`u $P$ est un polyn\^ome
holomorphe. Le fait que $J^*\cap H_{n-2}$ soit invariant
par $\Lambda$  entra\^{\i}ne que 
$$l_j(\gamma_n
z_{n-1},\overline\gamma_n\overline z_{n-1},
\gamma_n z_n+z_{n-1}, \overline
\gamma_n\overline z_n+\overline z_{n-1}) = d l_j (z_{n-1},
\overline z_{n-1},z_n,\overline z_n) +\Im P(z_{n-1},z_n).$$
\par
Soient $\alpha_j$, $\overline \alpha_j$, $\beta_j$ les coefficients de
$z_{n-1}\overline z_n$, $\overline z_{n-1} z_n$ et $|z_n|^2$ de
$l_j$. On supprime les
termes harmoniques et les termes $|z_n|^2$ de l'\'equation
pr\'ec\'edente. On obtient
$$\overline \gamma_n\beta_j z_{n-1}\overline z_n +
\gamma_n \beta_j\overline z_{n-1}z_n +
(\gamma_n\alpha_j+\overline\gamma_n\overline
\alpha_j +\beta_j)|z_{n-1}|^2=0$$
Par cons\'equent, $\beta_j=0$. Donc les $l_j$ sont harmoniques en $z_n$.
Posons $H_{n-1}:=\{z_1=\cdots =z_{n-1}=0\}$.
Alors $J^*\cap H_{n-1}$ est
d\'efini par $k-n$ \'equations du type $\Im z_j=\alpha_j
z_n^2+\overline \alpha_j \overline z_n^2$ et $J^*\cap H_{n-1}$
contient la courbe holomorphe d\'efinie par les \'equations
$z_j=2i\alpha_j z_n^2$.
Ceci contredit le fait que $J^*$ ne contient aucune
image holomorphe non constante 
de $\C$ et par cons\'equent $\Lambda_j$ ne d\'epend
pas de $z_s$ pour $s\not=j$ et $1\leq j\leq n$.
\par
Il reste \`a prouver que $\Lambda_j$ est ind\'ependante de $z_s$
pour tous $1\leq s\leq n$ et $n+1\leq j\leq k$.
On suppose qu'il existe $j\geq n+1$ et $1\leq s_1\leq s_2\leq n$
tels que $\Lambda_j$ contienne
le terme $z_{s_1}z_{s_2}$ avec un coefficient
$\alpha\not =0$. Comme $\Lambda$
est triangulaire, on a $\gamma_j=\gamma_{s_1}
\gamma_{s_2}$. La vari\'et\'e
$J^*$ est d\'efinie par $k-n$ \'equations $\Im z_s=
h_s(z_1,\overline z_1, \ldots, z_n,\overline z_n)$ pour
$j=n+1,\ldots, k$ o\`u les $h_s$ sont des polyn\^omes r\'eels
homog\`enes de degr\'e 2. En tenant compte l'invariance de $J^*$ 
par $\Lambda$ on a
$$h_j(\gamma_1 z_1, \overline \gamma_1 \overline z_1, \ldots,
\gamma_n z_n,\overline \gamma_n \overline z_n)=\Im \Lambda_j(\gamma_1
z_1,\ldots, \gamma_n z_n, h_j(z_1,\overline z_1,\ldots, z_n, \overline
z_n))$$
o\`u $\Lambda_j(z)= \gamma_j z_j+ P(z_1,\ldots,z_n)$ est ind\'ependant
de $z_{n+1}$, ..., $z_{j-1}$ et $P$ est un polyn\^ome.
Soit $\beta$ est le coefficient de $z_{s_1}z_{s_2}$
dans $h_j$. Les
coefficients de $z_{s_1}z_{s_2}$ des deux membres de
cette \'equation sont
$\gamma_{s_1}\gamma_{s_2}\beta$ et
$\gamma_j\beta+\gamma_{s_1}\gamma_{s_2}\alpha$. Or
$\gamma_j=\gamma_{s_1}\gamma_{s_2}\not =0$ et $\alpha\not=0$;
c'est la contradiction cherch\'ee.
\par
La vari\'et\'e $J^*$ est comme on a vu d\'efinie par $k-n$ \'equations $\Im
z_j=h_j(z_1,\overline z_1, \ldots, z_n,\overline z_n)$
o\`u les $h_j$ sont des
polyn\^omes r\'eels homog\`enes de degr\'e $2$ pour
tout $j$ avec $n+1\leq j\leq k$.
Posons $z':=(z_1,\ldots, z_n)$, $\overline z':=
(\overline z_1, \ldots,
\overline z_n)$ et $z'':=(z_{n+1},\ldots, z_k)$. Il existe des polyn\^omes
holomorphes homog\`enes
$P_j(z')$ de degr\'e 2 et des formes hermitiennes
$q_j(z',\overline z')$ tels
que $h_j=\Im P_j+q_j$. L'ensemble $J^*$ \'etant 
invariant par $\Lambda$, on a
$P_j(\gamma_1 z_1,\ldots, \gamma_n z_n)=\gamma_j
P_j(z_1,\ldots,z_n)$ et $q_j(\gamma_1 z_1,\ldots, \overline
\gamma_n \overline z_n)=\gamma_j q_j(z_1,\ldots,\overline z_n)$.
Quitte \`a effectuer le changement de coordonn\'ees $z_j\mapsto
z_j-P(z')$, on peut supposer que $P_j=0$ pour tout $j\geq n+1$. Comme $J^*$
ne contient aucune droite complexe, l'ensemble $\{z'|\ q(z',\overline
z')=0\}$ est \'egal \`a $\{0\}$ o\`u $q:=(q_{n+1},\ldots,q_k)$. 
Observons que si $f$ est
homog\`ene, le changement de coordonn\'ees ci-dessus 
ne modifie pas la rigidit\'e de $\varphi$.
\end{preuve}
%
%
%
%
%
%
%
\begin{proposition} Pour tout vecteur $c\in\R^{k-n}$,
la fonction $G^*$ est constante sur la vari\'et\'e
$J(c):=\{\Im z''=q(z',\overline z')+c\}$. Par cons\'equent, il existe
une fonction continue $\Phi$ telle que $\Phi(\Im z''-q(z',\overline
z'))=G^*(z)$.  
\end{proposition}
\begin{preuve} Par le choix de l'ouvert $\Omega$ contenant le point
  fixe $b$, la fonction $G^*$ et la
mesure $\mu^*$ sont simultan\'ement lamin\'ees au voisinage de $0$.
Consid\'erons les plaques de la lamination comme des graphes au-dessus
de $H$, l'espace tangent \`a $J^*$ en $0$. 
Il existe un voisinage $X$ de $0\in \C^n\times \R^{k-n}=H$,
un voisinage $Y$ de $0\in \R^{k-n}$, une application
r\'eelle analytique $\psi(z',\Re z'',v)$ de $X\times
Y$ dans $\R^k$ v\'erifiant les propri\'et\'es suivantes:
\begin{enumerate}
\item $\psi(z',\Re z'',0)=q(z',\overline z')$;
\item $\psi(0,0,v)=v$ pour tout $v\in Y$;
\item $G^*$ est r\'eelle analytique sur $X_v:=\{\Im z''=
\psi(z',\Re z'',v)\}$;
\item Pour tous $w_1$, $w_2$ dans $X_v$, on a
$G^*(w_1)\leq \delta(\|w_1-w_2\|)G^*(w_2)$ o\`u
$\delta$ est une fonction r\'eelle positive ind\'ependante de $v$
et tendant vers 1 en $0$.
\end{enumerate}
\par
L'application $\psi$ s'\'ecrit sous la forme:
$$\psi(z',\Re z'',v)=v+q(z',\overline z')+
v\O(\|z'\|)+v\O(\|\Re z''\|).$$
Posons $w:=(0,0,c)\in J(c)$. Soit
$x=(\alpha,\beta,\gamma)\in J(c)$. On a $\gamma=c+q(\alpha,
\overline \alpha)$. 
Il suffit de montrer que  $G^*(x)=G^*(w)$. Soit $s$ un
entier suffisamment grand. On pose $v:=cd^{-s}$,
$w_s:=(0,0,v)$, 
$(\alpha_s,\beta_s,0):=\Lambda^{-s}(\alpha,\beta,0)$,
$\gamma_s:=\psi(\alpha_s,\beta_s,v)$,
$y_s:=(\alpha_s,\beta_s,\gamma_s)$ et 
$x_s:=\Lambda^s(y_s)$. On a que $\Lambda^{-s} (J(c))=J(d^sc)$. 
En utilisant le d\'eveloppement
ci-dessus de $\psi$,
on montre facilement que
$x_s$ tend vers $x$ quand $s$ tend vers l'infini. Comme
$G^*$ est continue, il suffit de montrer que $\lim
G^*(x_s)=G^*(w)$. On a 
$G^*(x_s)=d^sG^*(y_s)$ et $G^*(w)=d^s
G^*(w_s)$. D'apr\`es la condition 4, on a
$$\delta(\|y_s-w_s\|)^{-1}G^*(x_s)\leq G^*(w)\leq
\delta(\|y_s-w_s\|) G^*(x_s).$$
Quand $s$ tend vers l'infini, on obtient 
$G^*(x)=G^*(w)$ car $y_s$ et $w_s$ tendent vers $0$ et
$\delta(\|y_s-w_s\|)$ tend vers 1.
\end{preuve}
\par
Nous voulons \`a pr\'esent construire les \'el\'ements du groupe $\A$
qui op\`ere transitivement sur les fibres de $\varphi$. Il est clair
que pour deux points $p$, $q$ d'une m\^eme fibre $G^*(p)=G^*(q)$ et
que pour des points $p$, $q$ g\'en\'eriques il existe une application
holomorphe $g$ telle que $g(p)=q$, $\varphi\circ g=\varphi$ et $G^*\circ
g=G^*$. On veut \'etudier les d\'eriv\'ees de $g$ \`a l'aide de la
relation pr\'ec\'edente. Cependant $G^*$ n'est pas d\'erivable, cela
oblige  \`a quelques d\'etours. On va \'etudier d'abord la restriction
de $G^*$ aux droites complexes.
\par
Soit ${\cal
D}\subset\{z'=0\}$ une droite complexe passant par $0$. On dit
que ${\cal D}$ est {\it non g\'en\'erique} si
$\dim_\R{\cal D}\cap\{z'=0,\Im z''=0\}=1$. 
Notons ${\cal D}^+\subset{\cal D}$ et
${\cal D}^-\subset{\cal D}$ deux
demi-plans complexes dont le bord est $\{z'=0, \Im
z''=0\}\cap {\cal D}$. 
\begin{corollaire} Pour toute droite non g\'en\'erique ${\cal
D}\subset\{z'=0\}$,
il existe des constantes non n\'egatives $c^+$ et $c^-$
d\'ependant contin\^ument de ${\cal D}^+$ et ${\cal D}^-$ telles
que $G^*(z)=c^+\|\Im z\|$ sur ${\cal D}^+$ et 
$G^*(z)=c^-\|\Im z\|$ sur ${\cal D}^-$. 
Pour toute droite complexe ${\cal
D}'\subset\{z''=0\}$ passant par $0$, il existe une constante
$c'>0$ d\'ependant contin\^ument de ${\cal D}'$ telle que
$G^*(z)=c'\|z\|^2$.
\end{corollaire}
\begin{preuve} Remarquons qu'une droite non g\'en\'erique rencontre
  $J^*$ le long d'une droite r\'eelle comme il r\'esulte des
  \'equations de $J^*$ ({\it voir} la proposition 5.6). 
\par
Observons que ${\cal D}$ est invariante par
$\Lambda$ car $\Lambda$ est diagonale. 
L'\'equation fonctionnelle $G^*(\Lambda z)=dG^*(z)$ et le lemme 5.5
  entra\^{\i}nent l'existence de $c^+$ et $c^-$.
Ces constantes d\'ependent contin\^ument de ${\cal
D}^+$ et ${\cal D}^-$ car $G^*$ est continue.
\par
Pour d\'eterminer $G^*$ sur une droite g\'en\'erique, on utilise la
connaissance de $G^*$ sur les droites non g\'en\'eriques et
l'invariance de $G^*$ sur $J(c)$. 
Fixons un $z=(z',0)\in{\cal D}'$. On pose $w:=(0,0,-q(z',\overline
z'))\in \{z'=0,\Re z''=0\}$. D'apr\`es la proposition 5.7, 
on a $G^*(z)=G^*(w)$.
Il existe un vecteur r\'eel $v\in \R^k$ qui ne d\'epend que de
${\cal D}'$ telle que $q(z',\overline
z')=v\|z\|^2$. D'apr\`es la partie pr\'ec\'edente, il existe une
constante $c>0$ telle que $G^*(z)=c\|v\|\|z\|^2$. Posons
$c':=c\|v\|$. On a $G^*(z)=c'\|z\|^2$.
Il est clair que $c'=c\|v\|$ d\'epend contin\^ument de ${\cal D}'$.  
\end{preuve}
\par
Nous allons \`a pr\'esent d\'efinir des ``d\'eriv\'ees'' de $G^*$ dans
certaines directions.
Soit $h$ une fonction d\'efinie sur une surface de Riemann $S$,
\`a valeurs dans $\C$.
Soit $z$ une coordonn\'ee locale de $S$ nulle en $a$. 
Si $u$ est un vecteur tangent en $a$ \`a $S$, 
il existe une constante $c$ telle que
$u=c\frac{\partial}{\partial z}$. Etant donn\'e une fonction
harmonique $l$, on
d\'efinit $u\otimes\overline u(h): =
c\lim_{z\rightarrow 0}[h(z)-l(z)]
|z|^{-2}$ lorsque cette limite existe. Observons que cette d\'efinition est
ind\'ependante de la coordonn\'ee $z$ et que lorsque la limite existe
elle est ind\'ependante $l$. C'est une formalisation de
$\partial^2/\partial u\partial \overline u$.   
\begin{lemme} {\bf i.} Soit $v$ un vecteur non nul de
$\{z'=0, \Re z''=0\}$. Soit $\sigma\subset
\C^k$ un arc r\'eel lisse issu de $0$ et tangent \`a $v$. 
Alors la d\'eriv\'ee
$v(G^*_{|\sigma})$ existe, elle est ind\'ependante de $\sigma$ et 
d\'epend contin\^ument de $v$. On la note $v(G^*)$. 
De plus
$v(G^*)=c^+ \|v\|$ o\`u $c^+$
est la constante associ\'ee \`a la demi-droite non g\'en\'erique
contenant $v$  qui est d\'efinie dans le corollaire 5.8. 
\par
{\bf ii.} Soit $u$ un vecteur holomorphe tangent 
\`a $\{z''=0\}$ en $0$.
Soit $S\subset \C^k$
une courbe holomorphe tangente \`a $u$ en $0$. Alors
$u\otimes \overline u(G^*_{|S})$ existe, elle est ind\'ependante de
$S$ et d\'epend contin\^ument de $u\otimes \overline u$. 
On la note $u\otimes \overline u(G^*)$
\end{lemme}
\begin{preuve} {\bf i.}
D'apr\`es la proposition 5.7, il suffit
de consid\'erer le cas o\`u $\sigma$ est contenue dans 
$\{z'=0, \Re z''=0\}$. Les \'equations de $J^*$ entra\^{\i}nent 
que $\{z'=0,\Im z''=0\}$
est contenu dans $J^*$. Par cons\'equent, $G^*$ s'annule
sur cet ensemble. Soient
$w=(w',w'')\in \sigma\setminus \{0\}$ et ${\cal D}_w$ la
droite complexe passant par $0$ et $w$. C'est une droite non
g\'en\'erique. On consid\`ere ${\cal
D}_w^+\subset {\cal D}_w$ le
demi-plan complexe contenant $w$ dont le bord est la droite
r\'eelle $\{\Im z''=0\}\cap{\cal D}_w$.
Notons \'egalement ${\cal D}$ (resp. ${\cal D}^+$)
la droite complexe (resp. le demi-plan complexe) qui contient $v$.
D'apr\`es le corollaire 5.8, il
existe une constante $c_w^+$
qui ne d\'epend que de ${\cal D}_w^+$ telle
que $G^*(w)=c_w^+|\Im w''|$. On a donc
$v(G^*_{|\sigma})=|v|\lim_{w\rightarrow 0}c_w^+$.
Comme $G^*$ est continue, cette limite existe et \'egale
\`a $|v|c^+$. Il est clair que cette constante ne
d\'epend que de $v$. De plus, elle d\'epend contin\^ument de $v$. 
\par
{\bf ii.} Pour simplifier les notations, on suppose que
$u=\partial/\partial z_1$. Alors la courbe $S$ est d\'efinie par
les \'equations $z_s=\psi_s(z_1)$ pour $s=2,\ldots,k$ o\`u
les fonctions $\psi_s$ sont holomorphes. Posons
$\psi':=(z_1,\psi_2,\ldots,\psi_n)$,
$\psi'':=(\psi_{n+1},\ldots,\psi_k)$ et $\phi:=(1,0,\ldots,0)$.
On note $\Pi$ la projection
$\Pi(z):=z_1$ et $G_S:=G^*\circ (\Pi_{|S})^{-1}$. 
D'apr\`es la proposition 5.7, on a
$G_S(z_1)=G^*(w)$ o\`u
$$w:=(w',\Re w'', \Im w'')=
(0,0,\Im\psi''(z)-q(\psi',\overline \psi')).$$
On note ${\cal D}^+$ la limite
de ${\cal D}_w^+$ quand $z_1\rightarrow 0$. Alors ${\cal D}^+$
est un demi-plan complexe qui ne d\'epend que de $u\otimes
\overline u$.
Sans perte de g\'en\'eralit\'e, on
peut supposer que ${\cal D}^+=\{z_1=\cdots=z_{k-1}=0, 
\Im z_k\geq 0\}$. Pour $z_1$ suffisamment petit,
$z_k$ est une coordonn\'ee de ${\cal D}_w$ et on a
${\cal D}_w^+=\{\Im z_k\geq 0\}\cap {\cal D}_w$. 
D'apr\`es le corollaire 5.8, il existe une constante $c_w$
d\'ependant contin\^ument de $w$ telle que 
$G_S(z_1)=c_w|\Im w_k|$ pour $z_1$ suffisamment petit.
On pose $c:=\lim c_w$.
\par
Comme $S$ est tangente \`a $u$, il existe des constantes $c_s$
telles que $\psi_s(z_1)=c_s z_1^2+\o(|z_1|^2)$ pour tout
$s=2,\ldots,k$. Posons $c'':=(c_{n+1},\ldots,c_k)$. Alors
$\Im w''=\Im c''z_1^2 -|z_1|^2q(\phi,\overline \phi)+\o(|z_1|^2)$.
D'o\`u
$$G_S(z_1)=c[\Im c_kz_1^2-|z_1|^2q_k(\phi,\overline
\phi)]+ \o(|z_1|^2).$$
On en d\'eduit que
$u\otimes \overline u (G^*_{|S})=
-cq_k(\phi,\overline \phi)$
car $\Im c_k z_1^2$ est harmonique. La continuit\'e de cette
quantit\'e est \'evidente.
\end{preuve}
\begin{lemme} Soit $K$ un sous-ensemble born\'e de $\R^n$
(resp. $\C^n$). Supposons 
$K$ non contenu dans aucun sous-espace r\'eel (resp. complexe)
propre de $\R^n$ (resp. $\C^n$).
Alors il existe un syst\`eme lin\'eaire de coordonn\'ees
dans lequel toute application lin\'eaire $L:\
\R^n\longrightarrow \R^n$ (resp. $L:\ \C^n\longrightarrow \C^n$)
v\'erifiant $L(K)=K$ est une isom\'etrie (resp. isom\'etrie
complexe).  
\end{lemme}
\begin{preuve} On consid\`ere le cas de $\C^n$; la preuve est
valable aussi pour le cas de $\R^n$. Soit $K'$ l'ensemble des
points $\lambda z$ o\`u $|\lambda|\leq 1$ et $z\in K$. Soit $H$
l'enveloppe convexe de $K'$. Comme $K$ est
born\'e et engendre $\C^n$,
$H$ est born\'e et d'int\'erieur non vide. 
Notons ${\cal O}_K$ (resp. ${\cal O}_H$)
le groupe des applications lin\'eaires complexes $L:\ \C^n
\longrightarrow \C^n$ v\'erifiant $L(K)=K$ (resp. $L(H)=H$).
On a ${\cal O}_K\subset{\cal O}_H$.
Comme $H$ est
born\'e et d'int\'erieur non vide,
le groupe ${\cal O}_H$ est compact. On en d\'eduit que ${\cal
O}_K$ est compact. On note
$e_1$, ..., $e_n$ les vecteurs de
la base orthonormale canonique de $\C^n$ et
$<z,z'>:=\sum_{i=1}^n z_i\overline z_i'$ le produit hermitien usuel
de $\C^n$. On d\'efinit 
$$<z,z'>_K:=\int_{{\cal O}_K}<L(z),L(z')> d\nu(L)$$ 
o\`u $\nu$ est la mesure de Haar de ${\cal O}_K$. On v\'erifie
facilement que $<.,.>_K$ est un produit scalaire et que
$<z,z'>_K=<L(z),L(z')>_K$ pour tout $L\in{\cal O}_K$.
Soit ${\cal M}$ la
matrice carr\'ee de rang $n$ d\'efinie par ${\cal M}:=(<e_i,e_j>_K)$.
Elle est d\'efinie positive. Il existe donc
une matrice inversible ${\cal N}$
v\'erifiant $\TF {\cal N}\overline {\cal N}={\cal M}$. Posons
$w:={\cal N}z$.
Pour ces nouvelles coordonn\'ees, on a
$$<w,w'>=<{\cal N}z,{\cal N}z'>=\TF z\ \TF {\cal N}\overline {\cal N} 
\overline z'
=\TF z {\cal M}\overline z'=<z,z'>_K.$$
Par cons\'equent, $<w,w'>=<L(w),L(w')>$ pour tout $L\in{\cal O}_K$. 
Dans ces nouvelles coordonn\'ees, $L$ est bien une isom\'etrie.  
\end{preuve}
Soit $N$ l'espace $\{z'=0, \Re z''=0\}$. 
Notons $K_N$ l'ensemble des vecteurs $v\in N$ v\'erifiant
$v(G^*)<1$ et $(-v)(G^*)<1$.
Soit $H$ l'espace complexe $\{z''=0\}$. On identifie $H$ avec l'espace
des vecteurs holomorphes tangents \`a $H$. 
Notons $K_H$ l'ensemble des vecteurs $u\in H$ v\'erifiant
$u\otimes\overline u(G^*)<1$. Les ensembles $K_N$ et $K_H$
repr\'esentent la variation de $G^*$ au voisinage de $0$ suivant
les directions tangentes de $N$ et de $H$.  
\begin{lemme} $K_N$ (resp. $K_H$)
est un ouvert born\'e de $N$ (resp. de $H$) 
contenant le point $0$.
\end{lemme}
\begin{preuve} 
D'apr\`es le lemme 5.9, $v(G^*)$
d\'epend contin\^ument de $v$. Par cons\'equent, $K_N$ est un
ouvert. Il est clair que $0\in K_N$. Il reste \`a montrer que
$K_N$ est born\'e.
\par
Sinon, soient $v_s\in K_N$ tels que
$\lim \|v_s\|=+\infty$. Notons ${\cal D}_s$ la droite complexe qui
est engendr\'ee par $v_s$. Notons \'egalement ${\cal D}_s^+$ (resp. ${\cal
D}_s^-$) le demi-plan complexe contenant $v_s$ (resp. $-v_s$)
dont le bord est $\{z'=0, \Im z''=0\}\cap {\cal D}_s$.
Alors pour tout $z\in {\cal D}_s^\pm$, on a 
$G^*(z)=c_s^\pm\|\Im z\|$ et $(\pm v_s)(G^*)=c^\pm_s\|v_s\|$. 
On en d\'eduit du lemme 5.9 que $\lim
c_s^+=\lim c_s^-=0$ car $\lim \|v_s\|=+\infty$ et
$(\pm v_s)(G^*)<1$.
La fonction $G^*$ \'etant continue, elle doit s'annuler sur toute
droite complexe adh\'erente \`a la suite $\{{\cal D}_s\}$.
C'est la contradiction cherch\'ee car $G^*$ ne
peut s'annuller sur aucune droite complexe.
\par
D'apr\`es le lemme 5.9, 
$u\otimes\overline u(G^*)$ d\'epend contin\^ument de $u$. Donc
$K_H$ est ouvert. Il est clair que $0\in K_H$. Il reste \`a
montrer que $K_H$ est born\'e.
\par
Sinon, d'apr\`es le corollaire 5.8,
pour toute droite complexe ${\cal D}\subset H$ il existe une
constante $c$ telle
que pour tout $z\in{\cal D}$ on ait $G^*(z)=c\|z\|^2$. Si
$K_H$ n'est pas born\'e, il existe des droites ${\cal D}_s\subset
H$ et des constantes $c_s$ tendant vers $0$ telles que
$G^*(z)=c_s\|z\|^2$ pour tout $z\in {\cal D}_s$.
Comme $G^*$ est continue, elle doit
s'annuller sur toute droite adh\'erente \`a la suite $\{{\cal D}_s\}$.
C'est la contradiction cherch\'ee.  
\end{preuve}
\par
D'apr\`es le lemme 5.11, quitte \`a effectuer des changements
lin\'eaires de
coordonn\'ees dans $\{z''=0\}$ et dans $\{z'=0\}$,
on peut supposer que toute
application lin\'eaire (resp. lin\'eaire holomorphe)
de $N$ (resp. $H$) pr\'eservant $K_N$
(resp. $K_H$) est une isom\'etrie (resp. isom\'etrie complexe).
Il est clair que
ces changements de coordonn\'ees pr\'eservent toutes les
propri\'et\'es de $J^*$, $G^*$ ainsi que la rigidit\'e de $\varphi$,
mais l'application lin\'eaire
$\Lambda$ n'est plus d\'efinie par une matrice diagonale. 
\par
Notons ${\cal M}$ la matrice du changement de coordonn\'ees
effectu\'e. Il existe des matrices carr\'ees $\M_1$ et $\M_2$ 
de rangs $n$ et $k-n$ telles que
$$\M=\left(
\begin{array}{cc}
\M_1 & 0\\
0   & \M_2
\end{array}
\right)$$
La matrice $\M\Lambda\M^{-1}$ est diagonale. Ses \'el\'ements
diagonaux sont \'egaux \`a $\pm d$ ou \`a $\sqrt{d}$ en module. On
constate que 
les $n$ premi\`eres (resp. $k-n$
derni\`eres) fonctions coordonn\'ees de $\Lambda$ sont
ind\'ependantes de $z'$ (resp. $z''$). Quitte \`a remplacer $\Lambda$
par $\Lambda^2$ on peut supposer que $\frac{1}{d}\Lambda''=\id$ o\`u
$\Lambda=(\Lambda',\Lambda'')$. D'apr\`es les propositions 5.7,
5.12 et le corollaire 5.8, 
l'application $\frac{1}{\sqrt{d}}\Lambda'$ est une isom\'etrie
complexe de $\C^n$. En effet,  $(\frac{1}{\sqrt{d}}\Lambda',\id)$
pr\'eserve la fonction $G^*$.
\par
Pour tout $u=(u',u'')\in J^*$,
on pose $\tau_u$ l'application lin\'eaire de
$\C^k$ dans lui-m\^eme d\'efinie par
$$\tau_u(z',z''):=(z'+u',z''+2iq(z', \overline u')+u'').$$
Alors $\tau_u(0)=u$.
D'apr\`es la proposition 5.7, on a $G^*\circ
\tau_u=G^*$ car $\tau_u$ pr\'eserve $J(c)$.
\begin{proposition}
Soit $W$ un ouvert de $\C^k$
rencontrant $J^*$. Soit $\tau=(\tau',\tau'')$
une application holomorphe ouverte de
$W$ dans $\C^k$ v\'erifiant $G^*\circ\tau=G^*$
dans $W$. Alors l'application
$\tau'$ ne d\'epend pas de $z''$ et elle
d\'efinit une isom\'etrie complexe de $\C^n$.
L'application $\tau$ est lin\'eaire. Il existe
$u\in J^*$ et $v\in J^*$ tels que $\tau_u\circ \tau$ et
$\tau\circ\tau_v$ soient 
lin\'eaires isom\'etriques. De plus, les $n$ premi\`eres (resp. $k-n$
derni\`eres) fonctions coordonn\'ees des applications  $\tau_u\circ \tau$ et
$\tau\circ\tau_v$ sont
ind\'ependantes de $z'$ (resp. $z''$). En particulier si $\tau(0)=0$
alors $u=v=0$ et $\tau'$, $\tau''$ sont lin\'eaires isom\'etriques.
\end{proposition}
\par 
Pour d\'emontrer cette proposition nous aurons besoin de quelques
pr\'eliminaires: 
\begin{lemme} Soit $\tau=(\tau',\tau'')$
une application holomorphe inversible d'un
voisinage de $0$ dans $\C^k$ v\'erifiant $\tau(0)=0$ et 
$G^*\circ\tau=G^*$.
Alors $\frac{\partial \tau'}{\partial z''}(0)=0$,
$\frac{\partial \tau''}{\partial z'}(0)=0$,
$\frac{\partial \tau'}{\partial z'}(0)$ et
$\frac{\partial \tau''}{\partial
z''}(0)$ sont
des isom\'etries complexes. De plus, $\frac{\partial
\tau''}{\partial z''}(0)$ est une matrice \`a coefficients r\'eels.
\end{lemme}
\begin{preuve} Comme $G^*\circ \tau=G^*$ dans un
voisinage $W$ de $0$, on a $\tau(J^*\cap W)=J^*\cap\tau(W)$. Donc
$\frac{\partial \tau}{\partial z}(0)$
pr\'eserve l'espace tangent complexe $H$ et l'espace tangent r\'eel
$L$ 
de $J^*$ en $0$. D'o\`u $\frac{\partial \tau'}{\partial z''}(0)=0$. La
d\'eriv\'ee $\frac{\partial \tau'}{\partial z'}(0)$
est une application
lin\'eaire de $H$. Cette application pr\'eserve $K_H$ car
$G^*\circ \tau=G^*$. C'est donc une
isom\'etrie complexe.
\par
Comme $\tau$ pr\'eserve $J^*$, on a $\Im
\tau''(z)=q(\tau'(z),\overline{\tau'(z)})$ pour tout
$z=(z',z'')\in W$ v\'erifiant $\Im z''=q(z',\overline z')$. 
Le d\'eveloppement de Taylor d'ordre 1 en $0$ des deux membres de
l'\'equation pr\'ec\'edente permet de voir que $\frac{\partial
\tau''}{\partial z'}(0)=0$. La d\'eriv\'ee $\frac{\partial
\tau''}{\partial z''}(0)$ 
pr\'eserve $\{z'=0, \Im z''=0\}$, elle
est donc une isom\'etrie complexe de $\C^{k-n}$ \`a
coefficients r\'eels et par suite elle  pr\'eserve $N:=\{z'=0, \Re z''=0\}$. 
Par cons\'equent,
elle pr\'eserve $K_N$ car $G^*=G^*\circ \tau$.
C'est donc une isom\'etrie 
complexe de $\C^{k-n}$ dont les coefficients sont r\'eels.   
\end{preuve}
{\it Preuve de la proposition 5.12---} 
Quitte \`a remplacer $W$ par un ouvert convenable,
on peut supposer que $\tau$ est injective.
Comme $G^*\circ\tau=G^*$, on a $\tau(J^*\cap
W)=J^*\cap\tau(W)$. 
Soient $w_1=(w_1',w_1'')\in J^*\cap W$ et
$w_2=(w_2',w_2''):=\tau(w_1)\in J^*\cap \tau(W)$.
Posons $\tilde\tau:=\tau_{w_2}^{-1}\circ \tau\circ\tau_{w_1}$.
On a $\tilde \tau(0)=0$, $\frac{\partial\tilde \tau'}{\partial
z}(0)=\frac{\partial\tau'}{\partial z}(w_1)$ et 
$G^*\circ\tilde \tau=G^*$. D'apr\`es le lemme 5.13,
$\frac{\partial\tilde \tau'}{\partial z''}(0)=0$.
D'o\`u $\frac{\partial
\tau'}{\partial z''}(w_1)=0$. Ceci est vrai pour tout $w_1\in
J^*\cap W$. C'est donc vrai pour tout $w_1\in W$ car $\tau$ est
holomorphe et $J^*\cap W$ est non pluripolaire. On en d\'eduit
que $\tau'$ est ind\'ependant de $z''$.
\par
D'apr\`es le lemme 5.13,
$\frac{\partial\tilde \tau'}{\partial z'}(0)$ est une isom\'etrie de $\C^n$.
Donc $\frac{\partial\tau'}{\partial z'}(w_1)$ est une isom\'etrie de $\C^n$
pour tout $w_1\in J^*\cap W$. Comme $\tau'$ est ind\'ependant de
$z''$, $\frac{\partial\tau'}{\partial z'}(w)$ est une isom\'etrie de $\C^n$
pour tout $w\in W$. Ceci implique que $\tau'$ est une isom\'etrie et
par suite $\tau'$ est affine.
\par
On va montrer que $\tilde \tau$ est une isom\'etrie. D'apr\`es la
partie pr\'ec\'edente, $\tilde\tau'$ est une isom\'etrie de
$\C^n$. D'apr\`es le lemme 5.13, il suffit de montrer que $\tilde
\tau''$ est une application lin\'eaire ind\'ependante de $z'$.
L'application $\tilde \tau$ pr\'eservant $J^*$, on a
$$\Im \tilde \tau''(z',\Re z''+iq(z',\overline
z'))=q(\tilde\tau'(z'),\overline{\tilde\tau'(z')}).$$
\par
Le membre de gauche est donc ind\'ependant de $z''$, 
$\tilde \tau$ \'etant holomorphe, $\tilde \tau''$ s'\'ecrit sous la forme
$\tilde \tau''(z)=\sigma(z')+l(z'')$ o\`u
$\sigma$ est une application holomorphe et $l$ est une
application lin\'eaire \`a coefficients r\'eels. On a donc
$$\Im\sigma(z')=-l(q(z',\overline
z'))+q(\tilde\tau'(z'),\overline{\tilde \tau'(z')}).$$
Par cons\'equent, $\Im \sigma(z')=0$ car le membre \`a droite ne
contient aucun terme harmonique. On a $\sigma=0$ et
donc $\tau''(z)=l(z'')$.
\par
Finalement, $\tau$ est lin\'eaire. Soient $u=\tau^{-1}(0)$ et
$v=\tau(0)$. Alors $\tau_u\circ\tau$ et $\tau\circ\tau_v$
v\'erifient les m\^emes propri\'et\'es que $\tilde\tau$. Ce
sont donc des applications lin\'eaires isom\'etriques. 
De plus, les $n$ premi\`eres (resp. $k-n$
derni\`eres) fonctions coordonn\'ees de ces applications sont
ind\'ependantes de $z'$ (resp. $z''$). 
\par \hfill $\square$
\begin{corollaire} Soit $W$ un ouvert de $\C^k$ rencontrant $J^*$.
Soit $\tau$ une application holomorphe ouverte de
$W$ dans $\C^k$ v\'erifiant $\varphi=\varphi\circ\tau$
dans $W$.
Alors $\tau$ est une application affine holomorphe.
De plus, l'ensemble
$\A$ de ces applications affines est un groupe et pour tout $\tau\in
\A$ on a $\tau(J^*)=(J^*)$.
\end{corollaire}
\begin{preuve} 
Comme 
$G^*=G\circ \varphi$, la relation $\varphi=\varphi\circ\tau$
implique $G^*\circ \tau=G^*$ sur $W$. D'apr\`es
la proposition 5.12, $\tau$ est une application affine
holomorphe.
\par
Par prolongement analytique, $\varphi=\varphi\circ\tau$ dans
$\C^k$. Comme $\tau$ est ouverte, elle est inversible. On a
$\varphi=\varphi\circ\tau^{-1}$. Donc $\tau^{-1}\in \A$. De
plus, si $\tau'$ est une application affine v\'erifiant
$\varphi=\varphi\circ\tau'$, on a
$\varphi=\varphi\circ\tau\circ\tau'$. Par cons\'equent, $\A$
est un groupe. La relation $J^*=\varphi(J)$ implique 
$$\tau(J^*)=\tau\circ\varphi^{-1}(J_k)=\varphi^{-1}(J_k)=J^*.$$
\end{preuve}
\par
En utilisant une id\'ee de Berteloot-Loeb \cite{BertelootLoeb}, on montre
la proposition suivante:
\begin{proposition}
Soit $\A_{|J^*}$ le groupe des automorphismes $\tau_{|J^*}$ de $J^*$
avec $\tau\in \A$. Alors $\A_{|J^*}$ est co-compact,
$\varphi(J^*)=J_k$ et $\A$ 
agit transitivement sur les fibres de $\varphi$. 
\end{proposition}
\begin{preuve} 
Par rapport \`a \cite{BertelootLoeb}, 
la seule diff\'erence dans notre cas est le fait que les \'el\'ements de
  $\A_{|J^*}$ ne sont pas tous isom\'etriques. La proposition 5.12 nous
  permet d'adapter l'id\'ee de Berteloot-Loeb. 
\par
L'application $\Lambda$ \'etant dilatante, 
pour montrer que $\A_{|J^*}$ est co-compact, il
  suffit de montrer qu'il existe un ouvert born\'e $W$ de $J^*$
  contenant $0$ tel que pour tout $w\in \partial W$,
il existe $w^*\in
  W$ v\'erifiant $\varphi(w^*)=\varphi(w)$. En effet, ceci implique
  que $W$ contient un domaine fondamental de $\A$.
\par  
Notons $\Pi$ la projection de $J^*$ dans $\C^n\times
  \R^{k-n}$ d\'efinie par $\Pi(z',z''):=(z',\Re z'')$.  
Soit $M>1$ tel que $\|q(z',\overline u')\|
\leq M\|z'\|\|u'\|$ pour tous $z'$
  et $u'$. On
choisit un nombre fini des points
  $a_s=(a_s',a_s'')\in J^*$ pour $s=1$, 2, ... v\'erifiant les
  propri\'et\'es suivantes:
\begin{enumerate}
\item Pour tout $s$, on a 
$1/2<\|\Pi(a_s)\|<1$.
\item Dans $\C^n\times \R^{k-n}$,
l'enveloppe convexe $U$ des points $\Pi(a_s)$ contient la boule de
centre $0$ est de rayon $1/2$. 
\item Pour tout $u\in \partial U$, il existe un $s$ tel que
  $\|u-\Pi(a_s)\|<1/16M$.  
\end{enumerate}
On choisit $b_0\in J^*$ tel que $\|b_0\|<1/100$ et tel que la mesure
$$\frac{1}{d^{sk}}\sum_{f^s(z)=\varphi(b_0)}\delta_z$$
tende vers $\mu$ quand $s$ tend vers l'infini. Alors l'ensemble
$\bigcup_{s\geq 0} f^{-s}(\varphi(b_0))$ est dense dans $J_k$. 
Quitte \`a pertuber l\'eg\`erement les points $a_s$, on peut supposer
qu'il existe un $p$ tel que $f^p(\varphi(a_s))=\varphi(b_0)$ pour tout
$s$. Posons $b_s:=\Lambda^p(a_s)$. Soit $V$ l'enveloppe convexe des
$\Pi(b_s)$ et $W:=\Pi^{-1}(V)$. 
\par
Pour $w\in \partial W$, soit $s$
tel que $\|\Pi(\Lambda^{-p}(w))-\Pi(a_s)\|<1/16M$. Si $\tau=(\tau',\tau'')\in
\A$ est tel que $\tau(b_0)=b_s$ et 
$w^*:=\tau^{-1}(w)$ alors $\varphi(w^*)=\varphi(w)$. On montre
que $w^*\in W$. 
Observons que $V$ est de taille $d^{p/2}$
dans les directions complexes et $d^p$ dans les directions r\'eelles.
D'apr\` es la condition 2, il suffit donc de montrer que
$\|{w^*}'\|<d^{p/2}/4$ et $\|\Re{w^*}''\|<d^p/4$. 
Comme $\tau'$ et $\frac{1}{\sqrt{d}}\Lambda'$ sont des isom\'etries de
$\C^n$ o\`u $\Lambda=(\Lambda',\Lambda'')$, 
on obtient 
$$\|{w^*}'-b_0'\|=\|w'-b_s'\|=d^{p/2}
\|\Lambda^{-p}(w)'-a_s'\|<d^{p/2}/16M\leq
d^{p/2}/16.$$
D'o\`u $\|{w^*}'\|<d^{p/2}/4$ car $\|b_0\|<1/100$.
\par
Posons $\tilde \tau=(\tilde
\tau',\tilde \tau''):=
\tau_{b_s}^{-1}\circ\tau\circ\tau_{b_0}$,
$w_1:=\tau_{b_s}^{-1}(w)$ et $w_2:=\tilde\tau^{-1}(w_1)$. 
On a alors $w^*=\tau_{b_0}(w_2)$. D'apr\`es la proposition 5.12
appliqu\'ee \`a $\tilde \tau$,
$\tilde \tau$ est lin\'eaire isom\'etrique et on obtient
$$\|w_2'\|=\|w_1'\|=\|w'-b_s'\|\leq
d^{p/2}/16M.$$
On a aussi $\Re w_1''=w''-\Re b_s''+2\Im q(w'-b_s',\overline
b_s')$. D'o\`u 
$$\|\Re w_1''\|\leq \|w''-\Re b_s''\|+
2M\|w'-b_s'\|\|b_s'\|\leq 3d^p/16$$
car
$$\|w''-\Re b_s''\|=d^p\|\Lambda^{-p}(w)''-a_s''\|<d^p/16M\leq
d^p/16$$
et
$$\|w'-b_s'\|\|b_s'\|=d^p \|\Lambda^{-p}(w)'-a_s'\|\|a_s'\|
\leq d^p/16M\leq d^p/16.$$
D'apr\`es la proposition 5.12 appliqu\'e \`a $\tilde \tau$, on a
$\|\Re w_2''\|\leq 3d^p/16$ car $w_2=\tilde \tau^{-1}(w_1)$.
Comme $w^*=\tau_{b_0}(w_2)$, on a $\Re {w^*}''=\Re
w_2''+\Re b_0-2\Im q(w_2',\overline b_0)$. Par cons\'equent,
$$\|\Re {w^*}''\|\leq \|\Re w_2''\|+\|\Re b_0''\| + 2M\|w_2'\|\|b_0\|\leq
\frac{3d^p}{16}+\frac{1}{100}+\frac{d^{p/2}}{800} \leq \frac{d^p}{4}.$$
\par
Il en r\'esulte que $\A_{|J^*}$ est co-compact. On a  
  $\varphi(J^*)=J_k$, en effet $\varphi(J^*)$ est un compact de $J_k$
  et contient l'ensemble dense $\bigcup_{s\geq 0}f^{-s}(\varphi(b_0))$. 
\par
Montrons que $\A$ agit transitivement sur les fibres de $\varphi$.
La proposition 5.12 montre que $\A$ agit
transitivement sur la fibre $\varphi^{-1}(z)$ pour tout $z\in J_k$. Comme $\#
  f^{-1}(z)\leq d^k$ pour tout $z$,
l'ensemble $\Lambda^{-1}(\A w)$ est une
r\'eunion d'au plus $d^k$ orbites de $\A$
pour tout  $w\in J^*$. Si $\varphi(w)\not\in f({\cal C}_f)$, on a
$\# f^{-1}(w)=d^k$ o\`u ${\cal C}$ signifie l'ensemble critique. La relation
$f\circ\varphi=\varphi\circ\Lambda$ implique que $\varphi^{-1}({\cal
C}_f)\subset \Lambda^{-1}({\cal C}_\varphi)$. 
Alors pour tout $w\in J^*\setminus
\Lambda^{-1}({\cal C}_\varphi)$, l'ensemble $\Lambda^{-1}(\A w)$ est une
r\'eunion de $d^k$ orbites de $\A$. 
Pour un $w\in J^*\setminus \Lambda^{-1}({\cal C}_\varphi)$ 
fix\'e, on peut choisir
  $\tau_1$, ..., $\tau_{d^n}$ des \'el\'ements de $\A$
  avec $\tau_1=\id$ 
  tels que $\bigcup\A.
  \Lambda^{-1}(\tau_s(w))=\Lambda^{-1}(\A. w)$.
\par
Par prolongement analytique, ceci est vrai pour tout $w\in \C^k$ \`a
l'exception d'une hypersurface complexe qu'on notera $S$.
Supposons
  que $\A$ n'agisse pas transitivement sur les fibres de
  $\varphi$. Il existe donc $w$ et $u$ tels que $u\not\in \A.w$
  et $\varphi(w)=\varphi(u)$. 
Quitte \`a pertuber l\'eg\`erement $w$ et $u$ on peut supposer que
$\Lambda^{-n}(\A w)\cap S=\emptyset$ et $\Lambda^{-n}(\A u)\cap S
=\emptyset$ pour
tout $n\geq 1$.
Par suite,
$\varphi(\Lambda^{-1}(\tau_s(w))$ (resp.
$\varphi(\Lambda^{-1}(\tau_s(u))$) sont les $d^k$ pr\'eimages
diff\'erentes de $\varphi(w)=\varphi(u)$ par $f$. On en d\'eduit
qu'il existe $s_1$ tel que
$\varphi(\Lambda^{-1}\circ\tau_1(w))=
\varphi(\Lambda^{-1}\circ\tau_{s_1}(w))$. Comme
$\tau_1=\id$, on a $\varphi(\Lambda^{-1}(w))=
\varphi(\Lambda^{-1}(\tau_{s_1}(w)))$.
\par
Par r\'ecurrence, il existe $s_1$, $s_2$, ... tels que pour tout
$n\geq 1$ on ait
$\varphi(w_n)=\varphi(u_n)$ o\`u $w_n:=\Lambda^{-n}(w)$ et
$u_n:=\Lambda^{-1}\circ
\tau_{s_n}\circ\cdots\circ\Lambda^{-1}\circ\tau_{s_1}(u)$. Comme
$u\not\in \A.w$, on a $u_n\not \in\A.w_n$ pour tout $n$.
La suite $u_n$ est born\'ee. En effet $\Lambda^{-1}$ est contractante
et la proposition 5.12 permet de
contr\^oler les $\tau_{s_j}$ qui sont proches d'isom\'etries (on
montre facilement que la composante en $z'$ est born\'ee; une
r\'ecurrence facile permet de v\'erifier que la composante en $z''$
l'est aussi). Soit
$u_0$ une valeur limite de cette suite. Comme $w_n$ tend vers
$0$,  $\varphi(u_0)=\varphi(0)$. Il existe donc $\tau\in \A$
tel que $\tau(u_0)=0$. La suite
$v_{n_r}:=\tau(w_{n_r})$
tend vers $0$ et v\'erifie
$\varphi(v_{n_r})=\varphi(w_{n_r})$. De plus puisque $u_n\not in \A
w_n$ on a $v_{n_r}\not=w_{n_r}$.
C'est la contradiction recherch\'ee car $\varphi$ est injective
au voisinage de $0$.
\end{preuve}
{\it Fin de la preuve du th\'eor\`eme 5.1---} Rappelons que nous devons
ici traiter le cas o\`u la p\'eriode $s$ de $b$ n'est pas \'egale \`a
1. Posons $F:=f^s$. Le point $b$ est fixe pour $F$. La fonction de
Green de $F$ est \'egale \`a celle de $f$. D'apr\`es les
lemmes et les propositions pr\'ec\'edents, on peut construire une
application holomorphe $\varphi$ et 
une application $\Lambda_F$ tels que $F\circ
\varphi=\varphi\circ\Lambda_F$ et tels que $\varphi(\C^k)$ soit de
compl\'ement pluripolaire. 
On peut \'egalement construire un groupe
$\A$ agissant transitivement sur les fibres de $\varphi$. Nous allons
montrer que $\varphi$ convient pour $f$.
\par 
On choisit $w_1$ et $w_2$ deux points de $J^*$
v\'erifiant $f\circ\varphi(w_1)=\varphi(w_2)$. Quitte \`a
pertuber l\'eg\`erement ces points, on peut supposer qu'ils
n'appartiennent pas \`a l'ensemble critique de $\varphi$. Soit
$\Lambda:= \varphi^{-1}\circ f\circ\varphi$ l'application
holomorphe 
d\'efinie au voisinage de $w_1$ v\'erifiant $\Lambda(w_1)=w_2$.
Alors $G^*\circ
\Lambda= dG^*$ car $G \circ f=d G$. Il suffit de
montrer que $\Lambda$ est affine. 
Posons
$\overline\Lambda:=\tau_{w_2}^{-1}\circ \Lambda\circ \tau_{w_1}$. 
\par
On a $G^*\circ
\overline \Lambda=dG^*$ et $\overline \Lambda(0)=0$. Il reste
\`a montrer que $\overline \Lambda$ est lin\'eaire. Posons $L(z',z''):=
(\sqrt{d}z',dz'')$. D'apr\`es le
corollaire 5.8, on a
$G^*\circ L=dG^*$. En effet, d'apr\`es la proposition 5.7, il suffit de
v\'erifier cette relation sur les droites non g\'en\'eriques.
\par
Posons
$\tau:=L^{-1}\circ \overline \Lambda$. On a
$G^*\circ\tau=G^*$. D'apr\`es la proposition 5.12, $\tau$ est
une isom\'etrie complexe. Donc $\overline \Lambda$ est lin\'eaire et
$\Lambda$ est affine. Par prolongement analytique, on a  
$f\circ\varphi=\varphi\circ \Lambda$ dans $\C^k$.
\par
\hfill $\square$
\begin{corollaire} Pour tout $z\in J_k$ l'ensemble
$\bigcup_{s\geq 0} f^{-s}(z)$ est dense dans $J_k$.
\end{corollaire}
\begin{preuve} Supposons que $\bigcup_{s\geq 0} f^{-s}(z)$ ne soit pas
dense dans $J_k$. Posons $U:=J_k\setminus\overline{\bigcup_{s\geq
0} f^{-s}(z)}$. Alors $U$ est un ouvert non vide. D'apr\`es la
proposition 5.15, il existe un ouvert non vide
$V\subset U$ qui est une vari\'et\'e r\'eelle analytique lisse.
Soit $b\in V$ un point p\'eriodique r\'epulsif de $f$.
D'apr\`es la proposition 5.15, toute voisinage de $b$
contient une pr\'eimage de $z$ par un it\'er\'e de $f$. Donc $U$ est vide.
\end{preuve}
{\it Fin de la preuve du th\'eor\`eme 1.1---}
La preuve du th\'eor\`eme 1.1 se termine exactement de m\^eme mani\`ere
que celle du th\'eor\`eme 5.1. Supposons d'abord que les $f_i$ soient
polynomiales. 
Dans les notations de la proposition 4.4, 
on choisit un
point $b\in J_k\cap \Omega$ 
p\'eriodique r\'epulsif pour $f_1$. Posons $F:=f_1^s$ o\`u
$s$ est la p\'eriode de $b$. On peut construire l'application
$\varphi$, l'application $\Lambda_F$ et le groupe $\A$ comme dans la
preuve du th\'eor\`eme 5.1. On peut \'egalement construire une application
affine $\Lambda_1$ v\'erifiant $f_1\circ \varphi=\varphi\circ
\Lambda_1$. Le fait que $f_1$ et $f_2$ ont la m\^eme fonction de Green
nous permet de construire $\Lambda_2$ de m\^eme
mani\`ere que celle de $\Lambda_1$.
\par
Le fait qu'on peut choisir $\varphi$ rigide dans le cas des applications
homog\`enes permet de passer du cas polynomial au cas g\'en\'eral.
\par \hfill $\square$
\section{Ensemble exceptionnel et orbifolds}
\ahead
Dans ce paragraphe, nous donnons quelques propri\'et\'es de l'ensemble
 exceptionnel puis nous d\'emontrons le 
 corollaire 1.3. Commen\c cons par la proposition
 suivante.
\begin{proposition} Sous l'hypoth\`ese du th\'eor\`eme 5.1,
 l'ensemble exceptionnel ${\cal E}:=\P^k\setminus \varphi(\C^k)$
 contient un nombre fini d'hypersurfaces complexes. La r\'eunion
 ${\cal E}'$ de telles hypersurfaces est totalement invariante,
 c'est-\`a-dire $f^{-1}({\cal E}')={\cal E}'$. Si $k=2$ l'ensemble ${\cal E}$
 est alg\'ebrique. 
\end{proposition}
\begin{preuve}
Notons $\E^*$ l'ensemble de $z\in \P^k$ tel que  la mesure
$$\mu^z_s:=\frac{1}{d^{ks}}\sum_{f^s(w)=z}\delta_w$$
ne tende pas vers la mesure d'\'equilibre $\mu$ quand $s\rightarrow
+\infty$. D'apr\`es \cite{BriendDuval2}, $\E^*$ contient un nombre
fini d'hypersurfaces alg\'ebriques. L'ensemble $\E$ \'etant totalement
invariant par $f$, le support $J_k$ de $\mu$ \'etant compact dans
$\P^k\setminus \E$ (proposition 5.15), $\E$ est donc contenu dans
$\E^*$. Par cons\'equent, il contient un nombre fini d'hypersurfaces
alg\'ebriques. La r\'eunion $\E'$ de ces ensembles alg\'ebriques est
totalement invariante par $f$ car $\E$ l'est.
\par 
Supposons maintenant que $k=2$.  D'apr\`es \cite{BriendDuval2}, 
${\cal E}^*$ est le plus
grand sous-ensemble analytique propre de $\P^2$ totalement invariant
par $f$. 
On a ${\cal E}\subset {\cal E}^*$. Il suffit de
montrer que ${\cal
E}^*={\cal E}$. Supposons que ${\cal E}^*\not={\cal E}$.
Avec les notations du
paragraphe pr\'ec\'edent, le point p\'eriodique r\'epulsif 
$b$ appartient \`a l'ensemble
$\overline{\bigcup_{s\geq 0} f^{-s}(c)}\subset{\cal E}^*$ pour tout
$c\in \E^*\setminus \E$.
D'apr\`es le corollaire
5.16, $\bigcup_{n\geq 0}f^{-n}(b)$ est dense dans $J_k$. Par
cons\'equent, $J_k$ est contenu dans
${\cal E}^*$. C'est la contradiction cherch\'ee car
$J_k$ n'est pas pluripolaire. On a donc ${\cal E}^*={\cal
E}$.  
\end{preuve}
{\it Preuve du corollaire 1.3---}
On construit maintenant l'orbifold ${\cal O}=(\P^k,n)$
associ\'e \`a $f_1$ et $f_2$. 
\par
Soit $H$ une hypersurface irr\'eductible de $\P^k$. On pose
$n(H)=\infty$ si $H\subset {\cal E}$ et $n(H)=m$ si $H\not\subset
{\cal E}$ et si la multiplicit\'e de $\varphi$ sur
$\varphi^{-1}(H)$ est \'egale \`a $m$. D'apr\`es le proposition 6.1, il
y a qu'un nombre fini d'hypersurfaces $H$ v\'erifiant $n(H)=\infty$. 
D'apr\`es le corollaire 5.14, si $H\not\subset {\cal E}$, 
les multiplicit\'es de $\varphi$ sur les composantes irr\'eductibles
de $\varphi^{-1}(H)$ sont \'egales car on passe de l'une \`a l'autre
par une application de $\A$. La fonction $n$ est donc 
bien d\'efinie. 
\par
Montrons maintenant que $f_i$ d\'efinit un
rev\^etement de ${\cal O}$ dans lui-m\^eme. Il est clair que 
$f_i$ d\'efinit un rev\^etement ramifi\'e de
$\P^k\setminus\bigcup_{n(H)=\infty} H$ dans lui-m\^eme. Soit
$H\not\subset {\cal E}$
une hypersurface irr\'eductible de $\P^k$. Il reste \`a montrer
que $n(f_i(H))=\mult(f_i,H)n(H)$. Posons $K:=\varphi^{-1}(H)$,
$L:=\Lambda_i(K)$.
D'apr\`es la relation $f_i\circ\varphi =\varphi\circ
\Lambda_i$, on a 
$f_i\circ\varphi(K)=\varphi\circ\Lambda_i(K)$. En comptant les
mutiplicit\'es, la derni\`ere relation nous donne
$$\mult(\varphi,K)\mult(f_i,H)=\mult(\varphi,L).$$
D'o\`u $n(H)\mult(f_i,H)=\mult(\varphi,L)=n(f_i(H))$ car
$L=\Lambda_i(K)=\Lambda_i\circ\varphi^{-1}(H)\subset \varphi^{-1}\circ
f_i(H)$. 
\par \hfill $\square$
\par
On obtient aussi avec la m\^eme d\'emonstration le corollaire suivant:
\begin{corollaire} Sous l'hypoth\`ese du th\'eor\`eme
5.1, il existe un orbifold ${\cal O}=(\P^k,n)$ tel que
$f$ d\'efinisse un rev\^etement de ${\cal O}$ dans lui-m\^eme.
\end{corollaire}
\section{Endomorphismes permutables de $\P^2$}
\ahead
Dans ce paragraphe, nous donnons la preuve du corollaire
1.11. L'id\'ee est que si l'ensemble exceptionnel $\E=\P^2\setminus
\varphi(\C^2)$ est vide le
groupe $\A$ est co-compact et le couple $(f_1,f_2)$ se trouve dans
l'exemple 1.10; sinon on montre que $\E$ contient une droite complexe
et on se ram\`ene au cas polynomial o\`u les endomorphismes
permutables sont d\'etermin\'es dans \cite{Dinh2}.
\par
Notons $F_1$, $F_2$ les relev\'es permutables de $f_1$, $f_2$ et
et $J_F$ l'ensemble Julia d'ordre $3$ de $F$. 
Cet ensemble est de dimension r\'eelle 5, 4
ou 3.
\par
Consid\'erons le cas o\`u $J_F$ est de dimension r\'eelle 5.
L'ensemble de Julia d'ordre 2 
de $f_1$ et $f_2$ est \'egal \`a $\pi(J_F)=\P^2$ o\`u
$\pi:\C^3\setminus\{0\}\longrightarrow \P^2$ est l'application
canonique. En effet, 
d'apr\`es la proposition 5.15, l'application $\varphi$
(d\'efinie pour $f_1$ et $f_2$) est surjective sur $J_2$, qui est
\'egal \`a $\P^2$ car 
il contient un ouvert. Le groupe d'applications affines $\A$ de
$\C^2$ qui agit transitivement sur les fibres de $\varphi$, 
est co-compact. On d\'eduit de la  
construction du paragraphe 5 
que les \'el\'ements de $\A$ sont des
isom\'etries complexes de $\C^2$ ({\it voir} la proposition 5.12).
Le couple $(f_1,f_2)$ se trouve
donc dans l'exemple 1.10.
\par
Supposons maintenant que $J_F$ soit de dimension 4. 
Montrons que
$\dim {\cal E}\geq 1$. 
D'apr\`es les paragraphes
pr\'ec\'edents, il existe une application lin\'eaire
diagonale $\Lambda$ et un nombre r\'eel $\alpha\not=0$
tels que $f\circ\varphi=\varphi\circ \Lambda$, 
$J^*:=\varphi^{-1}(J_2)=\{\Im
z_2=\alpha |z_1|^2\}$ o\`u $f$ est un certain it\'er\'e de $f_1$. Par
suite $\dim J_2=3$. Les
\'el\'ement de $\A$  du type $\tau\circ\tau_v$ o\`u $\tau$ est une
application lin\'eaire isom\'etrique complexe, $v=(v_1,v_2)\in J^*$ et 
$$\tau_v(z_1,z_2):=(z_1+v_1,z_2+2i\alpha \overline v_1 z_1+v_2).$$ 
Observons que $\tau_v$ pr\'eserve les vari\'et\'es $\{\Im
z_2=\alpha |z_1|^2+\const\}$.
%
%
\par
Quitte \`a effectuer un changement de coordonn\'ees du type
$(z_1,z_2)\mapsto (z_1,\pm z_2)$,
on peut supposer que $\alpha\geq 0$.
Posons $B_n:=\{\|z\|\leq n\}$ et  
${\cal D}_m:=\{z_2=-2mi\}$.
Notons ${\cal D}^*_m:=\A {\cal D}_m$.
Les \'el\'ements de $\A$ sont d\'ecrits ci-dessus. On montre
facilement que ${\cal D}^*_m\cap B_n=\emptyset$ pour $m>>n$.
D'apr\`es le corollaire 5.14 et la proposition 5.15,
$\varphi^{-1}\circ\varphi({\cal D}_m)= {\cal D}^*_m$. Par
cons\'equent, on a $\varphi(B_n)\cap \varphi({\cal D}_m)=\emptyset$ pour
$m>>n$ car $\varphi^{-1}(\varphi({\cal D}_m))\cap B_n=\A{\cal D}_m\cap B_n
=\emptyset$. 
Le compl\'ementaire de $\varphi(B_n)$ pour $n$ assez grand. On en
d\'eduit qu'il contient une image
holomorphe non constante 
de $\C$: $\varphi(\D_m)$. Donc $\dim {\cal E}\geq 1$.
De m\^eme mani\`ere, on montre que si $\dim J_F=3$ et $\dim J_2=2$ on a $\dim
{\cal E}\geq 1$.
\par
Notons ${\cal F}$ la courbe alg\'ebrique contenue dans
${\cal E}$, c'est-\`a-dire on ne consid\`ere pas les points isol\'es
\'eventuels 
de ${\cal E}$ ({\it voir} la proposition 6.1). 
Alors $f_1^{-1}({\cal F})=f_2^{-1}({\cal F})={\cal F}$.
D'apr\`es \cite{FornaessSibony1,CerveauLinsNeto}, 
${\cal F}$ contient une droite complexe.
\par
Si ${\cal F}$ contient trois droites complexes, $f_1$ est \'egale \`a
$[\lambda_0w_{\alpha_0}^{d_1}:\lambda_1w_{\alpha_1}^{d_1}:
\lambda_2w_{\alpha_2}^{d_1}]$
et
$f_2$ est \'egale \`a
$[\beta_0w_{\nu_0}^{d_2}:\beta_1w_{\nu_1}^{d_2}:
\beta_2w_{\nu_2}^{d_2}]$ 
o\`u $(\alpha_0,\alpha_1,\alpha_2)$,
$(\nu_0,\nu_1,\nu_2)$ sont des permutations de $(0,1,2)$ et
$\lambda_i$, $\beta_i$ sont des constantes non
nulles. On en d\'eduit que $(f_1,f_2)$ est conjugu\'e \`a un
couple dans l'exemple 1.6.
\par
Supposons que ${\cal F}$ contienne exactement deux droites complexes. 
Notons ${\cal E}_1$ et ${\cal E}_2$ ces deux
droites. Alors $f_i^2({\cal E}_j)={\cal E}_j$. Par
cons\'equent les $f_i^2$ sont polynomiales. D'apr\`es
\cite{Dinh2}, il existe un syst\`eme de coordonn\'ees affines tel
que l'une de deux conditions suivantes soit vraie:
\begin{enumerate}
\item $f_i^2=(\lambda_i z_1^{d_i^2},\pm \T_{d_i^2}(z_2))$;
\item
$f_i^2=(\lambda_i z_1^{d_i^2},\pm z_1^{d_i^2/2}
\T_{d_i^2}(z_2/\sqrt{z_1}))$
\end{enumerate}
o\`u les $\lambda_i$ sont des constantes non nulles.
\par
Il existe alors des constantes non nulles $\beta_i$
et des fonctions rationnelles $R_i$ telles que 
$f_i=(\beta_i z_1^{\pm d_i}, R_i(z_1,z_2))$. 
\par
Dans le premier cas,
comme l'application
$f_i^2$ pr\'eserve la fibration $\{z_1=\const\}$, elle pr\'eserve
aussi la fibration $f_i(\{z_1=\const\})$. On en d\'eduit que ces
deux fibrations sont \'egales. Par cons\'equent,
$R_i$ est ind\'ependante de $z_1$. Les fractions $R_i$
sont permutables et $R_i^2=\pm \T_{d_i^2}$. Donc $R_i=\pm
\T_{d_i}$ et $(f_1,f_2)$ se trouve dans l'exemple 1.7. 
\par
Dans le deuxi\`eme cas, on montre comme ci-dessus que $f_i$ pr\'eserve la
fibration $\{z_1=cz_2^2\}$. Par cons\'equent, il existe des
fonctions rationnelles d'une variable, permutables, $S_i$ telles que
$R_i(z_1,z_2)=z_1^{\pm d_1/2} S_i(z_2/\sqrt{z_1})$ et
$S_i^2=\pm\T_{d_i^2}$. On en d\'eduit que $S_i=\pm \T_{d_i}$ et
on v\'erifie facilement que $(f_1,f_2)$ est conjugu\'e \`a un
couple obtenu dans l'exemple 1.9 pour $h_i=z^{\pm d_i}$ 
\cite{Dinh2}.
\par
Si ${\cal F}$ contient une seule droite complexe, 
le couple $(f_1,f_2)$ est conjugu\'e \`a un couple
d'applications polynomiales. D'apr\`es \cite{Dinh2}, il est
conjugu\'e \`a l'un des couples d\'ecrits dans les exemples 1.5-1.9. 
Tien-Cuong Dinh \hfill Nessim Sibony\\
Math\'ematique - B\^atiment 425 \hfill Math\'ematique - B\^atiment
425\\
UMR 8628, Universit\'e Paris-Sud \hfill UMR 8628, Universit\'e Paris-Sud \\
91405 ORSAY Cedex (France) \hfill 91405 ORSAY Cedex (France)\\
TienCuong.Dinh@math.u-psud.fr \hfill 
Nessim.Sibony@math.u-psud.fr\\
\end{document}